\documentclass{amsart}[12pt]
 \usepackage{geometry}        % See geometry.pdf to learn the layout options. There are lots.
 \usepackage{graphicx}
 \usepackage{amssymb}
 \usepackage{epstopdf}
\usepackage{amsmath}
\usepackage{amsfonts}
\usepackage{mathtools}
\usepackage{amsthm}
\usepackage[unicode,pdfusetitle]{hyperref}
\usepackage{enumerate}
\usepackage{indentfirst}
\usepackage{bm}
\usepackage{scalerel}

\usepackage{comment}

\newtheorem{thm}{Theorem}

\newtheorem{ex}[thm]{Example}
\newtheorem{lem}[thm]{Lemma}
\newtheorem{defn}[thm]{Definition}
\newtheorem*{defn*}{Definition}
\newtheorem{cor}[thm]{Corollary}
\newtheorem*{thm*}{Theorem}
\newtheorem*{prop*}{Proposition}
\newtheorem{theoremA}{Theorem}

\newtheorem{rmk}[thm]{Remark}
\numberwithin{thm}{section}

\newcommand{\calR}{\mathcal{R}}
\newcommand{\calQ}{\mathcal{Q}}
\newcommand{\calA}{\mathcal{A}}
\newcommand{\calB}{\mathcal{B}}

\newcommand{\calD}{\mathcal{D}}
\newcommand{\calL}{\mathcal{L}}

\newcommand{\calF}{\mathcal{F}}
\newcommand{\calN}{\mathcal{N}}

\newcommand{\N}{\mathbb{N}}
\newcommand{\R}{\mathbb{R}}

\newcommand{\Q}{\mathbb{Q}}
\newcommand{\acl}{\operatorname{acl}}
\newcommand{\im}{\operatorname{im}}
\newcommand{\interior}{\operatorname{int}}
\newcommand{\gr}{\operatorname{gr}}
\newcommand{\vPhi}{\varphi}

\newcommand{\cal}[1]{\ensuremath{\mathcal{#1}}}

\newcommand{\forkindep}[1][]{%
  \mathrel{
    \mathop{
      \vcenter{
        \hbox{\oalign{\noalign{\kern-.3ex}\hfil$\vert$\hfil\cr
              \noalign{\kern-.7ex}
              $\smile$\cr\noalign{\kern-.3ex}}}
      }
    }\displaylimits_{#1}
  }
}
\newcommand{\zerofork}[1][]{%
  \mathrel{
    \mathop{
      \vcenter{
        \hbox{\oalign{\noalign{\kern-.3ex}\hfil$\vert$\rlap{$^{0}$}\hfil\cr
              \noalign{\kern-.7ex}
              $\smile$\cr\noalign{\kern-.3ex}}}
      }
    }\displaylimits_{#1}
  }
}

\newcommand{\dom}{\operatorname{dom}}
\usepackage{color}

\newcommand{\dcl}{\operatorname{dcl}}

\newcommand{\tp}{\operatorname{tp}}
\newcommand{\eX}{\operatorname{exp}}

\newcommand{\calM}{\mathcal{M}}

\newcommand{\calG}{\mathcal{G}}

\newcommand\Ofork{\scaleobj{.7}{\zerofork}}

\newcommand{\noopsort}[1]{}

\providecommand{\MR}{\relax\ifhmode\unskip\space\fi MR }
\begin{document}

\title[Companionability Characterization]{Companionability Characterization for the Expansion of an O-minimal Theory by a Dense Subgroup}
\author[Block Gorman]{Alexi Block Gorman}
\email{blockgoa@mcmaster.ca}
\address{Department of Mathematics and Statistics, McMaster University}
\address{1280 Main Street West, Hamilton, ON L8S 4K1, Canada}
\date{\today}
\subjclass[2010]{Primary 03C10 Secondary 03C64}
%\keywords{o-minimality, model companion, generic reduct, interpolative fusion, NIP, NTP2}
\maketitle

\begin{abstract}
This paper provides a full characterization for when the expansion of a complete o-minimal theory, one that extends the theory of ordered divisible abelian groups, by a unary predicate that picks out a divisible, dense and codense group has a model companion.
This result is motivated by criteria and questions introduced in the recent works \cite{KTW18} and \cite{E18} concerning the existence of model companions, as well as preservation results for some neostability properties when passing to the model companion.
Examples are included both in which the predicate is an additive subgroup of a real ordered vector space, and where it is a multiplicative subgroup of the nonzero elements of an o-minimal expansion of a real closed field.
The paper concludes with a brief discussion of neostability properties and examples that illustrate the lack of preservation (from the base o-minimal theory to the model companion of the expansion we define)  for properties such as strong, NIP, and NTP$_2$, though there are also examples for which some or all three of those properties are preserved.\\
\smallskip
\noindent \textbf{Keywords.} o-minimality, model companion, generic reduct, interpolative fusion, NIP, NTP2
\end{abstract}

%----------------------------
\section{Introduction}
%----------------------------

The results in this paper follow in the spirit of ``Generic structures and simple theories,'' a seminal paper of Chatzidakis and Pillay.
In \cite{CP98}, they show that the extension of a theory $T$ with uniform finiteness
in the language $\calL$ to an $\calL \cup \{P\}$-theory, where $P$ is a unary predicate with no induced $\calL$ structure, has a model companion $T_P$.
Moreover, if $T$ is simple in the sense of Shelah \cite{S80}, then so is $T_P$.
Many works in this vein have recently explored questions related to generic predicates and properties preserved  when passing to a model companion, including \cite{E18}, \cite{H04}, \cite{KR18}, and \cite{KTW18}.
For definitions and an in-depth introduction to existential closedness, model completeness, and model companions, see \cite{CK90}.

This paper concerns the extensions of a complete o-minimal theory, one that extends the theory of abelian groups, to a language with a unary predicate that picks out a dense, codense, and divisible subgroup of a definable and connected unary open set on which there is a group structure induced by $\calL$.
The main result of this paper is a full characterization for when these extensions have a model companion.
This result is motivated by the criteria and questions concerning the existence of model companions introduced in \cite{KTW18} and \cite{E18}.
The significance of fully characterizing the existence of a model companion is that although the conditions for companionability in \cite{E18} and \cite{KTW18} have impressive breadth, the authors of both papers found it elusive to obtain comparably broad negative results.
A sufficiently robust negative result in Section \ref{(R,G)} yields the complete characterization.

Throughout this paper, ``definable'' means definable with parameters unless specified otherwise.
Let $\calM$ be an o-minimal expansion of a real closed field. 
We define an \emph{exponential function on $\calM$} as a definable differentiable function 
$f: M \to M_{>0}$ such that $f'=f$.
 A \emph{partial exponential function on $\calM$} is by definition a definable differentiable function 
$f: (-\epsilon, \epsilon)\to M_{>0}$ for some $\epsilon>0$ in $M$, such that $f'=f$. 

Let $T$ be a model complete and o-minimal theory extending the theory of real closed fields.
Let $T^{\times}_{\calG}$ be the extension of $T$ to the language $\calL_{\calG} = \calL \cup \{ \calG \}$, along with axioms stating that $\calG$ picks out a dense, codense, and divisible multiplicative subgroup of the positive elements.
Then the characterization is the following:

\begin{theoremA}[Theorem \ref{mult dichotomy}]\label{mult}
Let $\calM$ be any model of $T$.
Then the theory $T^{\times}_{\calG}$ as described above has a model companion if and only if there is no partial exponential function on $\calM$.
\end{theoremA}

The main result of this paper is more general; we relax the assumption that $\calG$ picks out a multiplicative subgroup and allow the predicate $\calG$ to pick out a subgroup with respect to any binary operation $*$ in the language $\calL$ of $T$ on a definably connected, $\calL(\emptyset)$-definable open subset ${\bf D}$ of the underlying set, provided that $({\bf D}, <, *)$ is an ordered divisible abelian group.
In practice, this means that ${\bf D}$ is an open interval or ray. 
For example, in the context of $T^{\times}_{\calG}$ above, ${\bf D} = \R_{>0}$ and $*$ is multiplication, which only has a group structure on $\R^{\times}$, not all of $\R$.

For $\calM \models T$ and $\epsilon>0$ in $M$, we say that a function $f:I \to M$ is a local $(\calG,*)$-endomorphism if $I$ is an interval of radius $\epsilon$ about the $*$-identity, and for all $x,y \in I$ with $x*y \in I$, we have that $f(x*y) = f(x)*f(y)$.
The following theorem was conjectured by Erik Walsberg, whose work with Tran and Kruckman in \cite{KTW18} motivated this investigation of when an o-minimal expansion of a group with a predicate for a suitable subgroup is companionable. 

\begin{theoremA}[Theorem \ref{gen dichotomy}] \label{EWconjecture}
Let $T$ be a complete o-minimal theory extending the theory of ordered abelian groups, with language $\calL$.
Let $T_{\calG}$ be the $\calL \cup \{ \calG \}$-theory extending $T$ that also stipulates the following:
\begin{enumerate}
\item{$\calG$ is a dense and codense subset of an $\calL(\emptyset)$-definable open and definably connected unary set ${\bf D}$,}
\item{$\calG$ and ${\bf D}$ are ordered divisible abelian groups with respect to $* \in \calL$.}
\end{enumerate}
Then $T_{\calG}$ has a model companion if and only if every family of local $(\calG,*)$-endomorphisms whose graphs are uniformly definable in a neighborhood of the $(\calG,*)$-identity has finitely many germs at that point.
\end{theoremA}

For the negative direction, the proof of Theorem \ref{EWconjecture} involves showing that the class of existentially closed models of $T_{\calG}$ is not an elementary class if there is a family of $\calG$-endomorphisms whose graphs are uniformly definable in a neighborhood of the group identity, and the $\calG$-endomorphisms have infinitely many germs at that point.
We obtain the positive direction by demonstrating that the abscence of any such uniformly definable family of $\calG$-endomorphisms suffices to satisfy a slight modification of the hypotheses of the model companion axiomatization in \cite{E18}.

There are many theories that one might naturally consider which fall on each side of the companionability dividing line. For example, the following examples have a model companion by the result above:
\begin{itemize}
\item the theory of  $(\R,<, 0, +, (kx)_{k \in K})$, where $K \subseteq \R$ is a subfield, expanded by a predicate $\calG$ for a dense, divisible proper subgroup of  $(\R,+,0)$;
\item the theory of $(\R,<, 0,1, +, \cdot, (x^k)_{k \in \R})$ expanded by a predicate $\calG$ for a dense, divisible proper subgroup of $(\R_{>0}, \cdot,1)$.
\end{itemize}
By contrast this is untrue for the following examples of the negative result:
\begin{itemize}
\item the theory of the expansion of $(\R, <,0,1, +, \cdot |_{[0,1]})$ by a predicate $\calG$ for a dense, divisible proper subgroup of  $(\R,+,0)$;
\item the theory of $(\R, <,0,1,+,\cdot,e^x)$ expanded by a predicate $\calG$ for a dense, divisible proper subgroup of $(\R_{>0}, \cdot,1)$.
\end{itemize}

Part of the motivation of both \cite{E18} and \cite{KTW18} was to give preservation and non-preservation results for neostability properties, including simplicity, NIP, and NTP$_2$.
For facts about and definitions of NIP and NTP$_2$ see \cite{S15}, and for an introduction to simplicity see \cite{S80}.
Similarly, results in this paper illustrate a lack of preservation for such neostability properties, while also giving examples where preservation does occur.
\begin{theoremA}[Theorem \ref{TP2}]\label{multTP2}
For $T$ the theory of a real closed field, the theory $T_{\calG}^{\times}$, where $\calG$ is a subgroup of the multiplicative group, has a model companion $T^*_{\calG}$, and $T^*_{\calG}$ has TP$_2$.
\end{theoremA}
The above result mirrors the behavior of fusions observed in \cite{KTW18}.
On the other hand, the following preservation results for strong and NIP (respectively) do hold.
For an introduction to strong theories and finite burden, including definitions and relevant facts, see \cite{A07}.

\begin{theoremA}[Theorem \ref{VSNIP}, Corollary \ref{VSn}, Remark \ref{VSinf}]\label{strongVS}
Let $K$ be an ordered subfield of $\R$ and let $VS_K$ be the theory of ordered vector spaces over $K$. 
\begin{itemize}
\item{If $K$ has finite dimension $n$ as a vector space over $\Q$, then $(VS_K)_{\calG}$ has a model companion $(VS_K)^*_{\calG}$, and $(VS_K)^*_{\calG}$ is strong and has burden $n+2$.}

\item{If $K$ has infinite dimension as a vector space over $\Q$, then $(VS_K)_{\calG}$ has a model companion $(VS_K)^*_{\calG}$, and $(VS_K)^*_{\calG}$ has NIP but is not strong.}
\end{itemize}
\end{theoremA}

%%%%%%%%% DO NOT DELETE!!!!%%%%%%%%%%%%%%%%%%
The structure of this paper is as follows.
The remainder of this section gives background, notation, and preliminaries.
In section \ref{(R,G)}, we establish Theorem \ref{EWconjecture} by describing the sufficient condition for the non-existence of a model companion, then prove both the sufficient condition for non-existence of a model companion and the sufficient condition for existence of a model companion as independent lemmas.
In section \ref{ex}, we prove Theorem \ref{mult} and offer examples of o-minimal theories for which expanding by one type of group results in no model companion, whereas expanding by another type of group does yield a model companion.
Finally, in section \ref{stability} we show that for $T$ an o-minimal theory that expands a real closed field, even if the theory $T_{\calG}$ has a model companion, that model companion has $TP_2$. 
Conversely, if $T$ is an o-minimal theory with only a vector space structure, then the model companion of the same kind of expansion has NIP, but depending on the base field it may or may not be strong, i.e. have finite burden.
%%%%%%%%%%%%%%%%%%%%%%%%%%%%%%%%%%%%%%

\subsection{Background} 
In \cite{E18}, inspired by the work of Chatzidakis and Pillay, Christian d'Elb\'ee gives a natural sufficient condition for when a geometric theory $T$ that is model-complete augmented by a predicate that picks out
the underlying set of a reduct of the theory (subject to certain conditions) has a model companion.
He also shows preservation of NSOP$_1$ under certain assumptions. 
For a definition and relevant facts about the property NSOP$_1$, see \cite{S96}.
D'Elb\'{e}e frames both the criterion for companionability and the preservation results in the most natural terms that still allow for a large modicum of generality.
He includes some negative results in \cite{E18}, but a complete characterization of when a ``generic reduct'' has a model companion remained intractable.

In \cite{KTW18}, authors Kruckman, Tran, and Walsberg generalize the work of d'Elb\'{e}e by introducing a new set up known as ``interpolative fusions'' which they use to produce even broader conditions for when certain kinds of first-order theories have a model companion.
In particular, they define the ``fusion'' of two theories under the assumption that the original theories exhibit certain compatibility properties.
Under the assumptions that they establish, the fusion of the two theories is the model companion for the union of their theories (over some ``base theory'' that lies in their intersection).
Indeed, the authors show in \cite{KTW18} that the existence of a fusion of two model-complete theories is equivalent to the existence of a model companion for their union.
In the case that one considers a pair of disjoint theories with uniform finiteness, the existence of a model companion is due to Peter Winkler in his thesis \cite{W75}.

The dichotomy we establish in Theorem \ref{EWconjecture} bears some resemblance to the dividing line known as ``polynomial boundedness'' from Miller in \cite{M94}, where he establishes the ``growth dichotomy.''
The growth dichotomy states that for any o-minimal expansion $\calM$ of a real closed field, either there exists an exponential function on $\calM$, or every definable function that is not eventually zero is asymptotically equivalent to a power function.
Miller and Starchenko establish a similar growth dichotomy for ``linearly bounded'' expansions of the reals as an additive group in \cite{MS98}.
However, a notable difference from Theorem \ref{EWconjecture}, and the reason that neither growth dichotomy can be utilized directly to obtain the results in this paper, is that the companionability dividing line coincides with local definability of exponentiation or multiplication, rather than global definability or boundedness.

\subsection*{Notation.}
For brevity, we define $[n] = \{1, \ldots ,n\}$.
Given a tuple $\vec{x} = (x_1, \ldots ,x_n)$, we write $\hat{x}_i = (x_1, \ldots , x_{i-1}, x_{i+1}, \ldots ,x_n )$ to denote the vector of length $|\vec{x}|-1$ that is $\vec{x}$ with the $i^{th}$ component removed.
Additionally, if $\calM \models T$ and $Y \subseteq M^m$ is a definable set and $I \subseteq [m]$, define $\pi ^{I} (Y)$ to be the projection of elements in $Y$ onto the coordinates that are in the subset $I$.
When the indexing set $I$ is the singleton $\{i\}$, we will use $\pi^i$ rather than $\pi^{\{i\}}$.
Below we let $B^n_{\epsilon}(x)$ denote $\prod_{i=1}^n (x- \frac{\epsilon}{2}, x+ \frac{\epsilon}{2})$, i.e. the product of $n$ intervals of width $\epsilon$ centered at $x$ in $M^n$.
We let $\overline{U}$ denote the topological closure of the set $U \subseteq M^n$.
For a definable function $f$, we define the ``delta function,'' written $\Delta_t f(x)$, as $f(x+t) - f(x)$ for $t$ in a neighborhood of $0$ for which $x+t$ is within the domain of $f$.
We say that $f$ is definable or $\calL(\emptyset)$-definable interchangeably, and this means that the graph of $f$ is a $\calL(\emptyset)$-definable set in $\calM$.

\subsection{Preliminaries}
Throughout this paper, let $T$ be a complete, model complete theory that is o-minimal and extends the theory of ordered divisible abelian groups, and let $\calL$ be the language of $T$.
Henceforth we also let $\calM \models T$.
We will often use without stating it explicitly that because the structure $\calM$ is o-minimal, it has the uniform finiteness property, i.e. eliminates ``$\exists^{\infty}$''. 
For $a \in M$, we let $B_{\epsilon}(a)$ denote a ball of radius epsilon centered at point $a$.
Throughout, when we write that an interval $(a,b)$ or $[a,b]$ is ``nontrivial,'' we mean that the interval has nonempty interior.
For the remainder of this paper, we use ``$\calL(\emptyset)$-definable'' and ``definable'' to mean the same thing.
If a set $X$ is definable with parameters $A \subseteq M$, we write ``$\calL(A)$-definable'' and ``$A$-definable'' interchangeably.
Whenever a set $X$ is $\calL_{\calG}(A)$-definable, we will refer to it as such and never simply call such a set ``definable'' or ``$A$-definable.''

Recall the following definition from Theorem \ref{EWconjecture}:
\begin{defn}
Let $\calG$ be a unary predicate, and let $\calL_{\calG} := \calL \cup \{ \calG \}$.
We define the theory $T_{\calG}$ as follows:
\begin{enumerate}
\item{$T \subseteq T_{\calG}$,}
\item{$\calG$ is a dense and codense subset of ${\bf D}$, an $\calL(\emptyset)$-definable open and definably connected unary set,}
\item{$\calG$ and ${\bf D}$ are ordered divisible abelian groups with respect to $* \in \calL$.}
\end{enumerate}
\end{defn}

We use $(\calM, \calG)$ to denote models of $T_{\calG}$.
The set ${\bf D}$ remains fixed, but ${\bf D}$ does not factor into the notation, because none of the properties of $T_{\calG}$ relevant to this paper depend on ${\bf D}$ explicitly.
We will refer to the $<$-interval topology on $\calM$ merely as its topology, and similarly we call the topology induced by $<$ on $\calG$ simply the topology on $\calG$.
For any function $f: M^n \to M^m$ we denote the graph of $f$ by $\gr (f) \subseteq M^{n+m}$, and we let $\operatorname{dom}(f)$ denote the domain of $f$.
For functions $f,g:M \to M$ we use $f \circ g$ to denote $x \mapsto f(g(x))$, and call this function the precomposition of $g$ with $f$.

\begin{defn} \label{germs} Let $a \in M$. Let $\mathcal{A} = (f_i)_{i \in I}$ be a uniformly definable family of unary functions such that $a \in \bigcap_{i \in I} \dom (f_i)$. Then for each $f^* \in \mathcal{A}$ the set
\begin{gather*}
G_{f^*} := \{ (\epsilon, \gr(f|_{B_{\epsilon}(a)})\subseteq M^3: f \in \mathcal{A} \land \exists 0< \delta < \epsilon \forall x \in B_{\delta}(a))(f(x) = f^*(x)) \}
\end{gather*}
is $\calL(a)$-definable.

Call $\mathfrak{G}:= \{G_f: f \in \mathcal{A} \}$ the family of germs of $\mathcal{A}$ at $a$.
\end{defn}

We say that $\tilde{f}$ and $f^*$ have equivalent germs at $a \in M$ if there exists $r>0$ such that for all $\epsilon<r$ and for all $f \in \calA$, $(\epsilon, \gr (f|_{B_{\epsilon}(a)})) \in G_{\tilde{f}}$ if and only if $(\epsilon, \gr (f|_{B_{\epsilon}(a)})) \in G_{f^*}$.
We will denote the equivalence class of the germs in $G \in \mathfrak{G}$ by the partial function $\hat{f}_{G}: B_{\epsilon}(a) \to M$, where $\epsilon >0$ is the supremum over all radii on which some element of $G$ is defined, and $f \in \calA$ is such that $G_{f} = G$ for some $a \in M$.
Since we can generally assume $a=0$ in our setting, we will often make no reference to that parameter.
Below, let $D$ be on open and connected definable subset of $M$.

\begin{defn}\label{EIC}
As before, fix $\calM \models T$.
\begin{enumerate}
\item{
Suppose that $f:D \subseteq M \to M$ is a definable function and $U \subseteq D$ is a definable neighborhood of $0 \in M$, and that for all $x,y \in U$ such that $x+y \in U$ we have $f(x+y) = f(x) + f(y)$.
Then we say that $f$ is a definable endomorphism on $U$, or local endomorphism if $U$ is not specified.}

\item{
Let $m,n \in \N$ be such that $m>1$. 
We call a definable function $h:D \subseteq M^{m} \to M^n$ a definable transformation endomorphic in coordinates, or definable EIC transformation for brevity, if for each $i \in \{1, \ldots ,m\}$: 
\begin{align*}\emph{(I)} \text{  }
\forall t >0 \forall \hat{a}_i, \hat{b}_i \in \pi^{[m] \setminus i}(D) \forall x \in \pi^i(D) & \Big [ (a_1, \ldots ,x, \ldots a_{m-1})  \in D \land (b_1, \ldots ,x,\ldots b_{m-1}) \in D\\ 
  \rightarrow & \Delta_t h(a_1, \ldots ,x , \ldots a_{m-1}) = \Delta_t h(b_1, \ldots ,x , \ldots b_{m-1}) \Big ]
\end{align*}
where $\Delta_{t} h(a_1, \ldots ,x , \ldots a_{m-1}) =  h(a_1, \ldots ,x+t, \ldots a_{m-1})- h(a_1, \ldots ,x , \ldots a_{m-1})$ (and similarly for $\hat{b}_i$), and also 
\begin{align*}\emph{(II)} \text{  }
\forall \hat{a}_i \in \pi^{[m]\setminus i}(D) \forall x,y \in \pi^{i}(D) & \Big [ (a_1, \ldots ,x, \ldots a_{m-1})  \in D \land (a_1, \ldots ,y, \ldots a_{m-1}) \in D \rightarrow \\ h(a_1, \ldots ,x -y, \ldots a_{m-1}) & = h(a_1, \ldots ,x, \ldots a_{m-1}) - h(0, \ldots ,y, \ldots ,0) \Big ].
\end{align*}}
\item{Call a definable set $X \subseteq M^n$ with $n>1$ a definable hyperplane if there is some definable EIC transformation $h:M^{n-1} \to M$ such that $X = \gr (h)$.}
\item{Let $A$ be the following $n \times m$ matrix:
\[ \begin{bmatrix}
    g_{1,1} & \ldots & g_{1,m}\\
    \vdots & \ddots & \vdots \\
    g_{n,1} & \ldots & g_{n,m}
  \end{bmatrix} \]
where $g_{i,j} (x)$ is a definable endomorphism on $\pi_j(D)$.
For $\vec{a} \in M^m$, set: 
$$A \vec{a}:=\left( \sum_{j=1}^m g_{1,j}(a_j), \ldots ,\sum_{j=1}^m g_{n,j}(a_j) \right).$$
Given a definable EIC transformation $h:D \subseteq M^m \to M^n$,
we say that $A$ is the matrix representation of $h$ if the coordinate functions of $A$ are given by the following:
$$g_{i,j} (t) = \pi^j( \Delta_t h(a_1, \ldots ,a_{i-1},y,a_i, \ldots ,a_{m-1})),$$
where $\vec{a} = (a_1, \ldots ,a_{m-1}) \in \pi^{[m]\setminus i}(D)$ can be chosen arbitrarily.
We see that $A \vec{x} = h(\vec{x})$ in this case.
}
\end{enumerate}
\end{defn}

What part (2) of the above definition says is that $h:M^m \to M^n$ is a ``definable EIC transformation'' if when you view it as a function only of the $i^{th}$ coordinate (with $i \leq m$) then even if we vary the values of the other coordinates, the result is an endomorphism and has the same ``behavior'' (as measured by the delta function).

We will formulate many of the results in terms of definable EIC transformations because quantifier elimination for ordered divisible abelian groups tells us that these are precisely the definable unary functions in the language $(+,0,1)$ up to affine shifts.
Hence, for any other $\calL$-definable unary function $f$, if we knew that $f(\calG) = \calG$ and $f(\calG^c) = \calG^c$ in every model of $T_{\calG}$, then this would have to follow from the axioms of $T_{\calG}$ itself.
Yet $T_{\calG}$ is axiomatized in such a way that for every definable function $f$ that is not definable in the group language, the property ``$f$ sends some element of $\calG$ to $\calG^c$'' is realizable.

To see why it makes sense that we define EIC transformations using the delta function, in the lemma below we illustrate the link between definable local endomorphisms and delta functions that are constant.

\begin{lem}\label{delta}
Suppose that $f:D \subseteq M \to M$ is definable and $R >0$ is such that $(-R,R) \subseteq D$ and $\{x+t : x \in D , t \in [0, R) \} \subseteq D$. 
For all $t \in [0,R)$, the function $\Delta_t f(x)$ is constant (with respect to $x$) on ${\bf D}$ and $f(0)=0$ if and only if $f$ is a definable endomorphism on $(-R,R)$.
\begin{proof}
For the forward implication, suppose that $f$ and $R>0$ satisfy the hypotheses of the lemma.
Let $g(t) := \Delta_t f(x) = f(x+t) - f(x)$, for which $x$ is chosen arbitrarily since $\Delta_t f$ does not depend on $x$. 
Suppose that $y_1,y_2 \in D$ are such that $0<y_1<R$ and $y_1+y_2 \in D$.
Then $f(y_2 +y_1) - f(y_2) = g(y_1) = f(0+y_1) - f(0) = f(y_1)$. 
This yields $f(y_1 +y_2) = f(y_1) + f(y_2)$ as desired.

For the backwards implication, suppose that $f$ is a definable endomorphism on $[0,R)$.
By definition of a local endomorphism $f(0)=0$.
Let $t \in (0,R)$, and let $x \in [0,R)$ be such that $x+t \in D$.
Then $\Delta_t f(x) = f(x+t) - f(x) = f(x+t-x)$ since $f$ is a local endomorphism, and $f(x+t-x) = f(t)$, which is constant with respect to $x$, as desired.
\end{proof}
\end{lem}

From the above lemma we can deduce that each definable EIC transformation can be represented by a matrix consisting of the coordinate-wise delta functions.
Furthermore, we will observe that the definable EIC transformations that send $\calG$ to itself in models $(\calM,\calG)$ of $T_{\calG}$ are precisely those whose delta functions are $\emptyset$-definable in $(\calG,0,+)$.

\begin{lem}\label{matrix}
Let $h:D \subseteq M^m \to M^n$ be a $\calL( \vec{c})$-definable EIC transformation.
\begin{enumerate}
\item{There is a unique matrix representation $A$ for $h$ (as defined in \ref{EIC}).}
\item{The functions $g_{i,j}(t)$ used to define the matrix $A$ are all definable in $(\calG, 0,+)$, i.e. are $\Q$-affine, if and only if $T_{\calG} \cup \tp (\vec{c}) \vdash h(\calG^m) \subseteq \calG^n \land h(M^m \setminus \calG^m) \subseteq  M^n \setminus \calG^n$.}
\end{enumerate}
\begin{proof}
For (1), given a definable EIC transformation $h$ we define the corresponding matrix $A$ by $A\vec{x} = (\sum_{j=1}^m g_{i,j}(x_j))_{i \in [n]}$ where $g_{i,j}(t)= \Delta_{t} \pi^j(h(a_1, \ldots ,x, \ldots a_{m-1}))$.
By definition of an EIC transformation, the coordinate functions of an EIC transformation are also EIC.
By conditions (I) and (II) of the definition of an EIC transformation, 
the functions $g_{i,j}(t)$ are well-defined functions of $t$ alone, i.e. $\hat{a}_i =(a_1, \ldots,a_{m-1}) \in \pi^{[m]\setminus i}(D)$ can be arbitrarily chosen, and $x$ is not a free variable in $g_{i,j}$.
From Lemma \ref{delta}, we deduce that $\Delta_t$ applied to the functions of the form $\pi^j(h(a_1, \ldots ,x, \ldots a_{m-1}))$ are equivalent to these same functions as local endomorphisms.
Hence $g_{i,j}(t)= \Delta_{t} \pi^j(h(a_1, \ldots ,x, \ldots a_{m-1})) = \pi^j(h(a_1, \ldots ,t, \ldots a_{m-1}))$.
So we can write $h$ as the sum of its coordinate functions in the following way: for $\vec{x} \in D$,
$$h(\vec{x}) = \left ( \sum_{j=1}^m g_{1,j}(x_j), \ldots ,\sum_{j=1}^m g_{n,j}(x_j) \right ).$$
This shows the desired relationship of $h(\vec{x})$ and $A\vec{x}$. 
The uniqueness of the matrix $A$ follows from $A$ being determined by $A\vec{e_i} = h(\vec{e_i})$ for any set of basis elements $e_1, \ldots ,e_m \in D^m$.

For (2), let $\calL_0 = (0,<,+)$.
Quantifier elimination for ordered divisible abelian groups tells us that $\calL_0$-definable unary functions are all of the form $f(x) = qx+g$ where $q \in \Q$ and $g \in \calG$.
Hence if each $g_{i,j}$ is definable in $(\calG,0,+)$, then because it is an endomorphism, it is of the form $x \mapsto qx$ with $q \in \Q$.
The forward implication is immediate from this.

For the other implication, suppose that $h:M^n \to M$ is $\calL(\vec{c})$-definable, but not $\calL_0(\vec{c})$-definable, and that $h(\vec{g}) \in \calG^n$ if and only if $\vec{g} \in \calG^m$.
Then for each $j \in [m]$ and $i \in [n]$, it follows that $g_{i,j}(x) \in \calG$ if and only if $x \in \calG$ for every $(\calM,\calG) \models T_{\calG}$.
Let us assume that for some $j \in [m]$ and $i \in [n]$ the function $g_{i,j}(t)$ is everywhere non-constant, 
and on every interval it is definable in $\calL(\vec{c})$ but not in $\calL_{0}(\vec{c})$.

For any $(\calM, \calG) \models T_{\calG}$, we can pass to an $|\calM|^{+}$-saturated elementary extension $\calN \succcurlyeq \calM$.
Let $B \subseteq N \setminus M$ be a maximal $\dcl_{\calL}$-independent set over $\calM$ that is both dense and codense in $D^{\calN}$. 
Note that a maximal $\dcl_{\calL}$-independent set over $\calM$ need not be dense and codense in some interval, but we can construct such a $B$ that is dense and codense with the use of iterated choice functions.
See \cite{BV16} or \cite{BGCH21} for similar constructions.
Define $\calG^{(\calN, \calG)} =( \calG^{(\calM, \calG)} \oplus_{b \in B} \Q b) \cap D^{\calN}$.
It follows from o-minimality that $(\calN, \calG^{(\calN,\calG)}) \models T_{\calG}$.
For $g_{i,j}(b)$ to be in $\calG^{(\calN,\calG)}$, it would have to be the case that $g_{i,j}(b)$ is in the $\Q$-linear span of $\calG^{(\calM,\calG)}\cup B$, by how we defined the interpretation of $\calG$ in $(\calN, \calG)$.
However, by $\dcl_{\calL}$-independence of $B$ over $\calG^{(\calM,\calG)}$, we know that $g_{i,j}(b)$ is not in the $\Q$-linear span of $\calG^{(\calM,\calG)}\cup B$ for any $b \in B$, unless it is in the $\Q$-linear span of $\calG^{(\calM,\calG)} \cup \{ b \}$.
Yet $g_{i,j}(b)$ can't be of the form $qb+g$ with $q\in \Q$ and $g \in \calG^{(\calM,\calG)}$, because $g_{i,j}$ is an endomorphism, so it would have to agree with $x\mapsto qx$ on an entire interval, contradicting that it is not $\calL_0(\vec{c})$-definable on any interval.
Hence $g_{i,j}(\calG^{(\calN,\calG)}) \not \subseteq \calG^{(\calN,\calG)}$, as desired.
\end{proof}
\end{lem}

The following can be viewed as a corollary to Lemma \ref{matrix}, and will prove useful in section \ref{badCase}.
\begin{cor} \label{endoFact}
Let $h:D \subseteq M^n \to M^m$ be a definable EIC transformation with $B^n_{\epsilon}(\vec{0}) \subseteq D$ for some $\epsilon>0$ and  $\vec{x} = (x_1, \ldots ,x_n) \in D$.
For $K\subseteq [n]$ with $|K|=k$, define $\vec{x}_K=(x_{K,1}, \ldots ,x_{K,n})$ by $x_{K,i}=x_i$ if $i\in K$ and $x_{K,i}=0$ if $i \not \in K$.
Let $K_1, \ldots K_n \subseteq [n]$ be such that $|K_i|=k<n$ for all $i \in [n]$ and $K_i \neq K_j$ for all $i \neq j$.
The value of $h(\vec{x})$ is uniquely determined by the following:
$$h(x_{K_1}), \ldots ,h(x_{K_n})$$
(or, equivalently, determined by $h(x_1,0,\ldots,0), \ldots ,h(0, \ldots ,0, x_n)$ uniquely).
\begin{proof}
For $i \in [n]$, let 
Since $h$ is endomorphic in each coordinate, we can define 
$$\hat{h}_i(\vec{x}) = h(x_1, \ldots,x_{n-1},0) - h(x_1, \ldots, x_{i-1},0,x_{i+1}, \ldots , x_n) = h(0, \ldots, 0,x_i,0, \ldots ,0,-x_n).$$
Consider $\tilde{h}( \vec{x}) = \sum_{i=1}^n \hat{h}_i = h(x_1, \ldots, \allowbreak x_{n-1},-nx_n)$. 
We can $\calL ( \emptyset )$-define $\tilde{h}(x_1, \ldots , x_{n-1}, -\frac{1}{n}x_n)$ on the same domain, and we observe that $\tilde{h}(x_1, \ldots , x_{n-1}, \frac{-1}{n}x_n) = h(x_1, \ldots , x_n)$ as desired.
\end{proof}
\end{cor}

%----------------------------------------------------------------------------------%
\section{Model Companion in the General Setting}\label{(R,G)}
%----------------------------------------------------------------------------------%

Recall that throughout this paper we assume that $T$ is an o-minimal theory that extends the theory of ordered divisible abelian groups (which we call ODAG) and is both complete and model complete in the language $\calL$.
We require that $T$ have an $\forall \exists$-axiomatization in $\calL$ because we rely on the fact that for such theories, which are sometimes called ``inductive,'' the existence of a model companion is equivalent to the class of existentially closed models of $T$ having a first-order axiomatization.
In the theory $T_{\calG}$, recall that we require that $\calG$ is divisible, hence $\calG \models \rm{ODAG}$ as well.
The assumption that $\calG$ is a divisible group is only needed in section \ref{badCase}, where we will utilize the $\Q$-vector space structure on $\calG$ to prove that the framework in \cite{E18} is applicable to the theory $T_{\calG}$ in certain cases.

\begin{rmk}[Fact 2.4, \cite{G05}]
Let $\calM \models T$ be a $|T|^+$-saturated model, and let $\vPhi$ be an $n+m$-ary $\calL$-formula.  Then the set of tuples $\bar{b} \in M^m$ for which there is $\bar{a} \in M^n$ with $ \models \vPhi ( \bar{a} ,\bar{b})$ and $a_i \not \in \dcl_{\calL} (\hat{a}_i \cup \bar{b})$ for each $i \in [n]$ is definable.
\end{rmk}

\subsection{The case that $\calM$ defines an infinite family of distinct germs of endomorphisms at zero.}
%%%*~*~*~*~*~*
We will first give a name to a property of a theory $T$ which we then prove precludes the theory $T_{\calG}$ from having a model companion.
Below, when we say ``endomorphism'' we mean an endomorphism with respect to the binary operation (i.e. the corresponding symbol in the language $\calL$) with respect to which $\calG$ is a subgroup of ${\bf D}$. 
We will use the symbol ``$+$'' for this binary operation and the language of additive groups throughout this section.

\begin{defn}\label{UEP}
We say that an o-minimal theory $T$ has \emph{UEP} (uniform endomorphisms property) if there is an $\calL$-formula $\vPhi(x,\vec{y},z)$ for which in every model $\calM \models T$, there is an infinite definable set $J \subseteq M^{|\vec{y}|}$ such that for each $\vec{c} \in J$ there exists $\epsilon > 0$ such that the formula $\vPhi(x,\vec{c},z)$ defines the graph of an endomorphism on a neighborhood of $0$ with radius at least $\epsilon$, and for no other $\vec{d} \in J$ does $\vPhi(x,\vec{d},z)$ have the same germ at zero as $\vPhi(x,\vec{c},z)$. 
\end{defn}

Observe that in practice the property UEP reduces to the case where $\vec{c}$ is a singleton, since we can use definable choice to define a path through the infinite set $J$ that is parameterized by a single interval. 
Hence without loss of generality we shall work only with definable families of functions that vary with respect to a single parameter, but the formulas which define their graphs may require additional, fixed parameters.
By taking an appropriate closed subset of $J$ if necessary, we can also assume that $J$ is topologically closed.

Suppose $T$ is as above, and that there is a $\emptyset$-definable family of partial functions $\calF_Y := \{ f_y: D \subseteq M \to M : y \in Y \}$ for which the definable family of germs $\mathfrak{G}$ of $\calF_Y$ is infinite.
Then by the uniform definability of $F_Y$ and by definable choice for $T$, we conclude that $\mathfrak{G}$ is also uniformly definable.
In fact, we can uniformly definably choose a representative partial function for each element of $\mathfrak{G}$. 
This gives us a definable sub-family $F_{Y'}$ of representative partial functions, i.e. each partial function $f_{y}$ does not have the same germ at zero as $f_{y'}$ for any $y' \neq y \in Y'$.
Hence the definition of UEP is equivalent to the same statement with the uniqueness requirement for parameter $\vec{c}$ replaced by the requirement that the family of germs $\mathfrak{G}$ contains infinitely many distinct equivalence classes.

\begin{lem}\label{neglem}
Suppose that $T$ has UEP.
Then there is a $\emptyset$-definable family $H$ of definable partial functions such that for some interval $I \ni 0$ and infinitely many $q \in \Q \cap (0,1)$ we have $x \mapsto qx|_{I} \in H$.
\begin{proof}
Suppose that UEP is witnessed in $\calM \models T$ by the infinite, definable family of local endomorphisms $F_Y := \{ f_y: D \subseteq M \to M : y \in Y \}$ where $Y\subseteq M$ is the parameter space for the family of partial functions.
We claim there exists some $a >0$ and some nontrivial interval $Y' \subseteq Y$ such that $[-a , a] \subseteq \bigcap_{y \in Y'} \dom f_y$.
Suppose not, i.e. that for every $\epsilon >0$ there are only finitely many $y \in Y$ such that $[-\epsilon, \epsilon] \subseteq \dom f_y$.
Then by uniform finiteness, there is some $N \in \N$ such that for every $\epsilon >0$, at most $N$ elements of $F_Y$ have the interval $[- \epsilon, \epsilon]$ contained in their domains.
Yet this implies there are at most $N$ functions in $F_Y$ with distinct germs at $0$, since the domain of every function $f_y \in F_Y$ contains an interval centered at zero in its domain.
Hence there must exist such an $a>0$ and such a $Y' \subseteq Y$, and without loss of generality we may take  $Y:=Y'$.

We now prove that o-minimality allows us to definably choose $i,s \in Y$ for which on the interval $[-a,a]$ the functions $f_s:[-a,a] \to M$ and $f_i:[-a,a] \to M$ are such that every point $\{(x,z): x \in [-a,a]\}$ between their graphs is in the graph of one of the other endomorphisms in $F_Y$.
Since UEP dictates that $F_{Y}$ has an infinite family of germs at $0$, we conclude that $\{f_y(x):y \in Y\}$ is infinite for all $x \in [-a,a]$. 
As infinite unary sets in o-minimal structures have definable nonempty interior, we deduce that for each $x \in [-a,a] \setminus \{0\}$, the set $\{f_y(x): y \in Y\}$ contains an interval.
We (definably) choose $f_s$ so that $f_s(a)$ is in the upper half of the right-most interval in $F_Y(a):=\{ f_y(a): y \in Y\}$, and $f_i$ so that $f_i(a)$ is in the lower half of the right-most interval of $F_Y(a)$.
We observe that if $i \neq s$ and there is $x_0 \neq 0$ such that $f_i(x_0) = f_s(x_0)=y_0$, then for every $q \in \Q \cap [-1,1]$ we have $f_i(qx_0)=qy_0=f_s(qx_0)$. 
So by o-minimality, the set $\{x:f_i(x)=f_s(x)\}$ contains an interval that includes zero.
Since there are infinitely many $y \in Y$ such that $f_i(a)<f_y(a)<f_s(a)$, a similar argument shows that all such $f_y$ cannot cross $f_i$ nor $f_s$ on $(0,a]$. 
Therefore $\{z: f_i(x)<z<f_s(x)\}$ is contained in $F_y(x)$ for all $x \in (0,a]$ by o-minimality.
Hence, o-minimal cell decomposition allows us to choose an $f_s$ and an $f_i$ for which $f_s(x)$ and $f_i(x)$ are in the interior of $F_Y(x)$ for all $x \in (0,a]$.

By UEP, we can find such functions $f_s$ and $f_i$ whose graphs do not coincide on a neighborhood of $0$.
By making $a$ smaller if necessary, we can ensure that the graphs of $f_i$ and $f_s$ on $[-a,a]$ only intersect at $0$.
By our choice of $i \neq s \in Y$ and $a$, we ensure that for all $x \in (0,a]$ we know  $f_i(x)<f_s(x)$.
Similarly we chose $f_i$ and $f_s$ so that $\{(x,z) :x \in (0,a] \land (f_i(x)<z<f_s(x)) \}$ is contained in $\{(x,z): x \in (0,a] \land \exists y \in Y (f_y(x) = z)\}$.
Note that by antisymmetry of endomorphisms (i.e. $f(-x)=-f(x)$), all properties we have argued for $f_i$ and $f_s$ on $(0,a]$ hold with signs reversed on $[-a,0)$.
Without loss of generality, we restrict to the case that the family of functions $F_Y$ are only those whose graphs lie between that of $f_i$ and $f_s$ on all of the interval $[-a,a]$.

We now define a new family of functions which will contain the endomorphism $x \mapsto qx$ for each $q \in \Q \cap (-1,1)$. Let $\delta = f_s(a) - f_i(a)$.
As in the above argument showing that by restricting $Y$ to some $Y'$ we can ensure that the domains of all $f_y$ contain some interval $[-a,a]$, we can also find a subinterval $\tilde{Y} \subseteq Y$ and some positive elements $\delta' \leq \delta$ such that for all $f_{\tilde{y}}$ with $\tilde{y} \in \tilde{Y}$, we know that $[-\delta', \delta '] \subseteq \im f_{\tilde{y}}$.
Without loss of generality we assume that $Y = \tilde{Y}$ and change $\delta$ so that now $\delta' < \max \{ f_s(x)-f_i(x) : x \in [-a,a] \} \leq \delta$.

For each $y \in Y$, let $g_y (x)= f_y(x) - f_i(x)$ and let $\tilde{g}_y(x) = -g_y(x)$.  It is clear that $\{g_y :y \in Y\}\cup \{\tilde{g}_y :y \in Y\}$ is also a definable family of endomorphisms on $[-a,a]$ with distinct germs at $0$.
Observe that the partial inverse $g_s^{-1}:[-\delta, \delta ] \to [-a,a]$ is a continuous endomorphism since $g_s$ is one. 
Hence the family $\{h_y:[-\delta,\delta] \to [-\delta,\delta] : y \in Y\} \cup \{\tilde{h}_y: [-\delta,\delta] \to [-\delta,\delta]  :y \in Y\}$ given by 
$$h_y(x) = g_y(g_s^{-1}(x)), ~  \tilde{h}_y(x) = -h_y(x)$$
is again a family of partial endomorphisms with distinct germs at $0$.
To see this, we observe that for any $x,z \in [-\delta,\delta]$ such that $|x-z|< \delta$, we have $\tilde{h}_y(x-z) = -g_y(g_s^{-1}(x-z)) = -g_y(g_s^{-1}(x)-g_s^{-1}(z))$ because $x,-z \in \dom g_s^{-1}$.
Since $g_s^{-1}(x),-g_s^{-1}(z) \in [-a,a] \subseteq \dom g_y$, we know
$-g_y(g_s^{-1}(x)-g_s^{-1}(z))= -g_y(g_s^{-1}(x)) + g_y(g_s^{-1}(z)) = \tilde{h}_y(x)-\tilde{h}_y(z)$, and similarly for $h_y(x-z)$.
We will denote this family $H = \{h_y : y \in Y \} \cup \{\tilde{h}_y :y \in Y\}$.

Finally, we now show that $H$ contains the map $x \mapsto qx$ on $[- \delta, \delta ]$ for all $q \in \Q \cap (-1,1)$.
Fix one $q \in \Q \cap (-1,1)$ and let $y_q \in Y$ be such that $h_{y_q}(\delta) = q  \delta$.
We note that such a $y_q$ must exist because we chose $a>0$ and $i,s \in Y$, such that
for all $x \in [-a,a]$ and for all $z \in [f_i(x),f_s(x)]$ there exists $y \in Y$ such that $f_y(x)=z$.
This ensures that same holds for all $x \in [-a,a]$ and all $z \in [0,g_s(x)]$, and similarly for all $x \in [-\delta, \delta]$ and all $z \in [\tilde{h}_{s}(x), h_s(x)] = [-\delta , \delta]$.
Letting $x= \delta$ and $z = q \delta$, we deduce the existence of such a $y_q \in Y$.

We observe that for any $r \in \Q \cap (0,1)$ we have $h_{y_q}(r \delta ) = r h_{y_q}(\delta) = qr \delta$, 
hence the definable maps $h_{y_q}$ and $x \mapsto qx$ agree on infinitely many points between $0$ and $\delta$.
Thus they must agree on an interval containing infinitely many points of the form $rq \delta$ with $r \in \Q \cap (0,1)$.
In particular we know that they agree on the interval $[\frac{\delta}{2}, \delta] \subseteq [0,\delta]$ by o-minimality.
By shifting the elements $h \in H$ to the left by $\frac{\delta}{2}$ and down by $h(\frac{\delta}{2})$, and by defining $h_{y}(-x) = -h_{y}(x)$, we ensure that $H$ contains the germ of $x \mapsto qx$ at zero, as desired.
We note also the $\emptyset$-definability of $H$ follows from the $\emptyset$-definability of $a$ and $\delta$, which is immediate by o-minimality and definable choice for $T$. 

\end{proof}
\end{lem}

With this lemma we can now prove the following negative result:

\begin{thm}\label{negprop}
Suppose that $T$ has UEP. 
Then there is no model companion for the theory $T_{\calG}$.

\begin{proof}
By Lemma \ref{neglem}, there is a definable family of partial functions $H = \{h_y:M \to M : y \in Y\}$ (where $Y$ is the parameter space for the family of functions) as described in the above lemma.
We now suppose for contradiction that the class of existentially closed models of $T_{\calG}$ is axiomatizable, say by the theory $T^{EC}_{\calG}$. 
Let $(\calM , \calG ) \models T^{EC}_{\calG}$ be an $\aleph_1$-saturated model.
We let $I \subseteq \bigcap\{ \dom(h_y) :y \in Y \}$ be nonempty and $\emptyset$-definable.
Note that such an interval $I \ni 0$ exists (and is nontrivial) by the remarks preceding Lemma \ref{neglem}.
We may further assume the interval $I$ is contained in ${\bf D}$.
Note that ${\bf D} = \interior (\overline{\calG})$ since $\calG$ dense and codense in ${\bf D}$ implies ${\bf D} \subseteq \interior (\overline{\calG})$, and ${\bf D}$ being an open and definably connected unary set in an o-minimal structure implies it must be an interval or ray at infinity (positive or negative), so $\overline{\calG}\setminus {\bf D}$ is the (at most two) endpoints of ${\bf D}$.

First, let us reindex the family $H$ using an element $\gamma \in \calG \cap I$ (by density of $\calG$ in $I$, can find such an element) in the following way.
We let $Y_{\gamma} = \{ h_y(\gamma) :y \in Y \}$, and by Lemma \ref{neglem} we note that for infinitely many $q \in \Q \cap (0,1)$ the element $q\gamma$ is in $Y_{\gamma}$.
Since we are simply reindexing, for each element $y^*:= h_y(\gamma) \in Y_{\gamma}$ where $y \in Y$, we still define $h_{y^*} (x) = h_y(x)$ for all $x \in \dom h_y$.
We let $X := \{ y \in Y_{\gamma} : \forall g \in ( \calG \cap I) ( h_y(g) \in \calG) \}$.
We define a superstructure $( \calM' , \calG') \supseteq (\calM, \calG)$ by taking $\calM \preccurlyeq \calM'$ to be an $|M|^+$-saturated elementary superstructure; let $I'$ be the interpretation of $I$ in $\calM'$, and let $\calG' = \calG \bigoplus_{b \in \calB } \Q b$ where $\calB \subseteq I'$ is  $\dcl_{\calL}$-independent over $\calM$.
Using choice, we may ensure that $\calB$ is both dense and codense in $I'$; see e.g. \cite{BGCH21} for such a choice construction.
By construction it is clear that $( \calM' ,\calG') \models T_{\calG}$, it is an extension of $(\calM, \calG)$ since $\calG'\cap \calM = \calG$,  and $(\calM, \calG)$ is existentially closed in $(\calM', \calG')$ since we assume $(\calM, \calG)$ is in the class of existentially closed models of $T_{\calG}$.

Let $X'$ be the interpretation in $(\calM ',\calG ')$ of the formula that defines $X$ in $(\calM,\calG )$.
Since for every $y \in X'$ we have $(\calM', \calG') \models h_y(b) \in \calG'$, we know $b \in \calB$ implies $h_y(b) = \sum_{i\in [m]} q_i b_i + a$ for some $b_1, \ldots ,b_m \in \calB$ and $q_1, \ldots ,q_m \in \Q$ and $a \in \calG$.
We suppose that $y \in X' \cap \calG '$, and note that there are infinitely many elements in $X' \cap \calG '$ since $\gamma \in X' \cap \calG '$, and so is every positive rational multiple smaller than $\gamma$.
Without loss of generality, we may assume that $h_y(x)$ is not identically zero on any nontrivial interval containing zero.

We know that $y = \sum_{i\in [m']} q'_i b'_i + a'$ for some $q'_1, \ldots ,q'_{m'} \in \Q$ and $b'_1, \ldots ,b'_{m'} \in \calB$ and $a' \in \calG$.
Lightly abusing notation, we will write $y=y(\vec{b}',a')$.
Hence for each $b \in \calB$ we conclude
$$h_{y(\vec{b}',a')}(b) = \sum_{i \in [m]} q_i b_i + a.$$
This in turn implies that $b \in \dcl_{\calL}(\{b_1, \ldots ,b_m,b'_1, \ldots ,b'_{m'}\} \cup \calG$), since $\dcl$ is a pregeometry in o-minimal structures.
This would contradict that $\calB$ is $\dcl_{\calL}$-independent over $\calG$ unless $h_y(b)$ is also $\calL( \{ b \} \cup \calG)$-definable.
Hence we conclude that $q_i \neq 0$ precisely if $b_i = b$.
Since $h_y(b) = q_ib + a$ we also know for each $n \in \N$ that $h_y(b/n) = q_ib/n + a/n$. 
By o-minimality the definable function $h_y(x) - q_ix$ is either non-constant with respect to $x$ or takes on finitely many values on any interval to the right of zero.
Since $a \in \calG$, it cannot be a non-constant function of $b \in \calB$, which means $a =0$.
Moreover, since $h_y$ agrees with $x \mapsto q_i x$ on all rational multiples of $b$, it agrees with this function on an interval by o-minimality.
Since this is true for each $b \in \calB$, we conclude that on a sufficiently small interval around zero (depending on $y$) we must have $h_y(x) = qx$ for some $q \in \Q$.
Since $y, y' \in Y$ and $y \neq y'$ implies $h_y \neq h_{y'}$, we conclude that $y = q \gamma$ for some $q \in \Q$.

Finally, we will observe that $X \cap \calG \subseteq X'\cap \calG$ as subsets of $M'$, and we conclude that $X \cap \calG$ is a countable set.
Since $\calM \preccurlyeq \calM'$, we have $Y_{\gamma}^{\calM} = M \cap Y_{\gamma}^{\calM'}$.
It only remains to see that $(\calM, \calG) \models \forall g \in \calG \cap I (h_y(g) \in \calG)$ implies that $(\calM', \calG') \models \forall g \in ( \calG \cap I)(h_y(g) \in \calG)$ for each $y \in Y \cap M$.
This is immediate from the existential closedness of $(\calM, \calG)$ in $(\calM', \calG')$, so $y \in X \cap \calG$ implies $y \in X'\cap \calG'$.
This contradicts the axiomatizability of the existentially closed models, since then $(\calM, \calG)$ must be an uncountably-saturated model which defines a countable infinite set.
\end{proof}
\end{thm}

\subsection{The case that definable families of endomorphisms in $\calM$ have finitely many germs.} \label{badCase}
%%%*~*~*~*~*

Throughout this section, let $T$ be a theory that does not have UEP. 
Let $\calM \models T$ be a saturated model.
To establish our criterion, we will examine the interaction of definable curves in $M^n$ with definable endomorphisms in multiple variables in the context of linearly bounded structures. 
Here, we use ``curve'' to mean a function $f: I \to M^n$ where $I \subseteq M$ is an open interval.

We are now ready to state and prove the demarcation lemma for companionability of $(\calM, \calG)$.
This lemma, and the theorem which follows, will prove Theorem \ref{EWconjecture} in conjunction with Theorem \ref{negprop} from the previous subsection.
In essence, the following lemma shows that we can recover any hyperplane that intersects an arbitrary definable curve in $M^n$ on an open subset of its domain.

\begin{lem}\label{demarcation}
Suppose $m,n \in \N$ and let $n \geq 2$. For every definable family of germs of curves at the origin:
$$\{F(x, \vec{a}) = (f_1(x, \vec{a}), \ldots ,f_n(x, \vec{a})): \vec{a} \in A \subseteq  M^m \}$$
there are only finitely many definable hyperplanes $H \subseteq M^{n}$ for which there exists an $n$-box $B^n_{\epsilon}(\vec{0})$ and a tuple $\vec{a} \in A$ such that $\calM \models \forall x \forall \vec{y} \in B^n_{\epsilon}(0) ( F(x, \vec{a})= \vec{y} \rightarrow \vec{y} \in H)$.
\begin{proof}
First, remark that we can extend the family of functions $\{F(x,\vec{a}):a \in A\}$ to the higher-dimensional family $\{F'(x,\vec{a}) = (x,F(x,\vec{a})): a \in A\}$ and in doing so regard the graphs of the original family of functions as the image of the new family of functions.
In light of this, we may address the case that some coordinate function $f_i(x,\vec{a})$ is itself a local endomorphism of $x$ by applying the statement of the lemma to the extension of $F(x,\vec{a})$ by the coordinate function $x \mapsto x$.
By the monotonicity theorem, we can define the finite set of intervals $\calD_{\vec{a}} \subseteq \dom F(x, \vec{a})$ on which each $f_i$ is injective and continuous, or constant. 
We will proceed by induction on $n$, the dimension of the image of $F$.
We exclude the intervals on which some $f_i$ is constant, since on those intervals the result will follow by induction hypothesis.

We perform a series of manipulations on $F(x, \vec{a})$, in which we iteratively precompose the inverse of the coordinate function of the $i^{th}$ coordinate and then take the $\Delta$-function of the resulting curve. 
We will do this for the first coordinate function, conclude that we can do the same for each coordinate, and the result will follow by induction.
We may further expand the family $F(x,\vec{a})$ and the parameter tuple $\vec{a}$ so that $F(x, \vec{a}b) := F(x - b, \vec{a}) -F(b,\vec{a})$ is in the family of functions for each $b \in \dom F(x,\vec{a})$ and each $\vec{a} \in A$.
Below we will suppress the parameter tuple $\vec{a}$, using $F(x) :=F(x, \vec{a})$, and $f_i(x) :=f_i(x, \vec{a})$.
For each $H$ a definable hyperplane that intersects the image of $F$ as specified in the hypotheses, let $h:M^{n-1} \to M$ be a definable EIC transformation witnessing that $H$ is such a definable hyperplane, as described in Definition \ref{EIC}.

If $n=2$, we can precompose $f_1^{-1}$ with $F(x)$ and conclude that on some interval $f_1(B_{\epsilon_a}(0))$ we can define the function $\tilde{F}(y) = F(f_1^{-1}(y))= (y, f_2(f_1^{-1}(y)))$.
If we assume that for some $\vec{a}$ the image of $F$ intersects the graph of a definable endomorphism $h:M \to M$ on some interval, then we conclude that 
$F(f_1^{-1}(y))= (y, f_2(f_1^{-1}(y)))= (y, h(y))$ on that interval.
By the definition of $T$ not having UEP, this $h$ can only be one of finitely many local definable endomorphisms for all $\vec{a} \in M^m$.
The case $n=2$ thereby reduces to the $n=1$ case, so we now take our base case to be $n=3$.

Consider 
$$F(x) = \left( f_1(x), f_2(x), h(f_1(x), f_2(x)) \right)$$ 
defined on $B_{\epsilon}(0)$, where $h:M^2 \to M$ is a definable EIC transformation.
As indicated by our choice of domain decomposition, we assume that $F$ is continuous and injective in each coordinate on $B_{\epsilon}(0)$.
We make the assumption that $f_i(0)=0$ for $i \in [2]$ since if some curve in the family $F$ nontrivially intersects a shift of a definable hyperplane, then the curve shifted to pass through the origin (which we have included in the family $F$) intersects the unshifted definable hyperplane.
We define $\tilde{F}(y) = (y,f_2 \circ f_1^{-1} (y), h(y, (f_2 \circ f_1^{-1}) (y)))$ on 
$B_{\epsilon_1}(0) \subseteq f_1(B_{\epsilon}(0))$ for a suitably chosen $\epsilon_1>0$.
We now consider the function $ \Delta_{t_1} \tilde{F} (y) = \tilde{F}(y+t_1) - \tilde{F}(y)$ and observe 
$$\tilde{F}(y+t_1) - \tilde{F}(y)= (t_1, f_2 \circ f_1^{-1} (y+t_1) - f_2 \circ f_1^{-1} (y), h(t_1, (f_2 \circ f_1^{-1}) (y+t_1)-(f_2 \circ f_1^{-1} )(y))).$$
Let $f_{2,-1}(x):= f_2(f_1^{-1}(x))$. 
There are now two cases, by the o-minimality of $\calM$.
Either for all $t_1>0$ sufficiently small there is some $\tilde{\epsilon}_1>0$ such that the function $\Delta_{t_1}  f_{2,-1} (y)$ is constant on $(0, \tilde{\epsilon}_1)$ with respect to $y$,
or there exists a $t_1>0$ such that for any sufficiently small interval $(0,\tilde{\epsilon}_1)$ the function $\Delta_{t_1} f_{2,-1} (y)$ is monotone increasing or monotone decreasing with respect to $y$ on $(0,\tilde{\epsilon}_1)$.

Suppose we are in the first case, i.e. for any sufficeintly small $t_1>0$ the function $\Delta_{t_1} f_{2,-1} (y)$ is constant on $(0, \tilde{\epsilon}_1)$ for some $\tilde{\epsilon}_1>0$.
We can choose $\tilde{\epsilon_1}$ to satisfy the hypothesis of Lemma \ref{delta}, and thus conclude that $\Delta_{t_1} f_{2,-1} (y)$ being constant on a neighborhood $(0,\tilde{\epsilon}_1)$ means that $f_{2,-1} (y)$ is itself a local endomorphism on $(0, \tilde{\epsilon}_1)$.
Consequently, we observe that $h(y , f_{2,-1}(y))$ intersects a definable hyperplane on $(0, \tilde{\epsilon}_1)$.
If there were infinitely many distinct choices (as a function of parameters $\vec{a}$) for $h$, then projecting $\tilde{F}(x)$ onto its first and third coordinates would yield a definable family of endomorphisms, contradicting the assumption that $T$ does not have UEP.
So we conclude that the claim holds for case 1.

Assume now we are in case 2. 
Without loss of generality, let $\tilde{\epsilon}_1>0$ be such that there exists $t_1 \in (0, \frac{\epsilon_1}{2})$ for which the function $\Delta_{t_1} \tilde{F} (y) = \tilde{F}(y+t_1) - \tilde{F}(y)$ is continuous and injective as a function of $y$ in each coordinate on interval $I_1 := (0,\tilde{\epsilon}_1)$.
Fix one such $t_1$, and we define this function as $F^{(1)} := \Delta_{t_1} \tilde{F} (y)$, which  equals
\begin{gather*} \tilde{F}(y+t_1) - \tilde{F}(y) = (y+t_1,f_{2,-1}(y+t_1),h(y+t_1,f_{2,-1}(y+t_1))) -(y,f_{2,-1}(y),h(y,f_{2,-1}(y))) \\
= (t_1,\Delta_{t_1} f_{2,-1}(y),h(t_1,\Delta_{t_1} f_{2,-1}(y))).
\end{gather*}
We remark that by our assumption, $\Delta_{t_1}f_{2,-1}(y)$ is invertible as function of $y$ on a neighborhood of $0$, so we now define $\tilde{F}' :\Delta_{t_1} f_{2,-1} (I_1) \to M^3$ given by \\
$$\tilde{F}'(x):= F^{(1)}( (\Delta_{t_1} f_{2,-1})^{-1} (x))= (t_1, x , h(t_1,x)).$$
Since we definably chose $t_1>0$, the tuple $(t_1, 0, h(t_1))$ is uniformly definable in the same parameters $\vec{a}$ as $\tilde{F}'$, hence the function $\hat{F}(x) = \tilde{F}'(x) - (t_1, 0, h(t_1,0))$ is uniformly definable in parameters $\vec{a}$.

Shifting $\hat{F}$ to have $0$ in its domain if necessary, we see that $\hat{F}$ projects onto the definable hyperplane given by $\tilde{h}_1(x) = (x, h(0,x))$, which is uniformly definable in terms of parameters $\vec{a}$. 
Exchanging the roles of $f_1$ and $f_2$ in the above proof, we similarly see that $\tilde{h}_2 = (x, h(x,0))$ is a hyperplane uniformly definable in parameters $\vec{a}$ as well.
By Lemma \ref{endoFact}, together $\tilde{h}_1$ and $\tilde{h}_2$ uniquely determine the definable EIC transformation $h$.
By induction hypothesis, there are only finitely many hyperplanes that $\tilde{h}_1(x)$ and $\tilde{h}_2(x)$ can agree with as $\vec{a}$ ranges over the parameter space, and hence there are finitely many definable EIC transformations that $h$ can possibly be for all parameter tuples $\vec{a}$.

The induction step proceeds analogously to the base case.
Let $n>3$ and suppose that the induction hypothesis holds for all curves in $M^n$.
Let $F(x) = (f_1(x), \ldots , f_n(x), h(f_1(x), \ldots, f_n(x)))$ be as in the statement of the lemma, with $h:M^n \to M$ a definable transformation.
We make the same assumption as in the base case, that on $B_{\epsilon}(0)$ each coordinate function is injective, and surjects onto a neighborhood of the origin of at least diameter $\epsilon$.
We define:
$$\tilde{F}_1(y) = \left( y,f_2 \circ f_1^{-1} (y), \ldots ,f_n \circ f_1^{-1}(y), h(y,f_2 \circ f_1^{-1} (y), \ldots ,f_n \circ f_1^{-1} (y) \right)$$
on $B_{\epsilon_1}(0) := f_1(B_{\epsilon}(0))$.
Note that $f_1$ is a continuous and injective map from $\calR$ to $\calR$, so it maps open balls to open balls. 
We can replace $f_1(x)$ with $f_1(x) - f_1(0)$, and let $\epsilon_1$ be half the diameter of $f_1(B_{\epsilon}(0))$.
Replacing $f_1^{-1}$ with $f_i^{-1}$, we can define a function $\tilde{F}_i$ analogously.

We now consider $F'(y) = \Delta_{t_1} \tilde{F}_1 (y) = \tilde{F}_1(y+t_1) - \tilde{F}_1(y) = (t_1, f_2 \circ f_1^{-1} (y+t_1) - f_2 \circ f_1^{-1} (y), \ldots ,f_n \circ f_1^{-1} (y+t_1) - f_n \circ f_1^{-1} (y), h(t_1,f_2 \circ f_1^{-1}(y+t_1)-f_2\circ f_1^{-1} (y), \ldots ,f_n \circ f_1^{-1} (y+t_1)-f_n \circ f_1^{-1} (y) )$.
Let $f_{i,-1}(y):= f_i(f_1^{-1}(y))$ for each $i \in \{2, \ldots ,n\}$.

We definably choose $t_1>0$ and $\epsilon_1>0$ such that $\Delta_{t_1} f_{i,-1}(y)$ is constant with respect to $y$ for the fewest possible number of indices $i \in \{2, \ldots ,n\}$ on $(0,\epsilon_1)$.
We define the components of a new tuple $\vec{g}(t_1)$ as follows:
$$g_i(t_1)= \begin{cases} 
t_1, & i=1 \\
\Delta_{t_1} f_{i,-1}(\frac{\epsilon_1}{2}), & \text{if } \Delta_{t_1} f_{i,-1}(y) \text{ is constant on an interval }(0, \delta), i \in \{2, \ldots ,n\} \\
 0,& \text{if } \Delta_{t_1} f_{i,-1}(y) \text{ non-constant on each interval }(0, \delta), i \in \{2, \ldots ,n\} \\
h(y, g_2(t_1), \ldots ,g_n(t_1)) & i=n+1 \end{cases}.$$
We observe that if $\Delta_{t_1} f_{i,-1}$ is not constant on any interval with left endpoint $0$ for any $2 \leq i \leq n$, then $g(\vec{t})= (t_1, 0, \ldots ,0, h(t_1,0, \ldots ,0))$.
We remark that the tuple $\vec{g} \in M^{n+1}$ is uniformly definable in the same parameters $\vec{a}$ as $F$.
Hence the function $\hat{F}(x) = F'(x) - \vec{g}(t_1)$ is uniformly definable over $\vec{a}$ as well.

In the case that $\vec{g(t_1)} = (t_1, 0, \ldots ,0, h(t_1,0, \ldots ,0))$ we observe the following:
$$\hat{F}(x) = \big ( 0, \Delta_{t_1} f_{2,-1} (y), \ldots , \Delta_{t_1} f_{n,-1} (y), h(0, \Delta_{t_1} f_{2,-1} (y), \ldots ,\Delta_{t_1} f_{n,-1} (y)) \big)$$ 
is equal to $\{0\} \times \hat{F}_0$ where $\hat{F}_0$ is a curve in $n-1$-dimensional space.
We apply the induction hypothesis to conclude that there are finitely many possible $n-1$-dimensional EIC transformations with which $\hat{F}_0$ can locally coincide.
Hence the $n^{th}$ coordinate of $\hat{F}$ can only coincide with one of finitely many EIC transformations with respect to the change of variables $\hat{z}_0:= (0,\Delta_{t_1} f_{2,-1} (y), \ldots , \Delta_{t_1} f_{n,-1} (y))$.
Let $\tilde{h}_1(\hat{z}_0):= h(0, z_1, \ldots ,z_{n-1})$ denote this coordinate function with the change of variables.

For each $f_k$ with $k \in \{2, \ldots ,n\}$, instead of precomposing $f_i$ with $f_1^{-1}$, we precompose $f_i$ with $f_k^{-1}$ and repeat the argument to conclude that there are similarly only finitely many possible $n-1$-dimensional EIC transformations with which $\hat{F}_k$, which is defined analogously to $\hat{F}_0$, can locally coincide.
Then we perform a similar change of variables to that in the paragraph above, and define $\tilde{h}_k( \hat{z}_k) = h(z_1, \ldots ,z_{k-1},0,z_{k+1},\ldots ,z_n)$ where $z_i := \Delta_{t_k} f_{i,-k} (y)$ for each $i \in [n] \setminus \{k\}$.
By Lemma \ref{endoFact}, an EIC transformation is uniquely determined by how it acts on each coordinate, so we can recover the EIC transformation $h$ uniquely from the coordinate functions $\tilde{h}_1, \ldots ,\tilde{h}_n$.
We conclude that there are only finitely many $n$-dimensional EIC transformations that $h$ can be.

Now consider the case that $g_{k}(t_1) \neq 0$ for the indices $k \in \{2, \ldots ,n\}$ such that $k \in J$, where $\emptyset \neq J \subseteq \{2, \ldots n\}$.
Recall that this means $\Delta_{t_1}f_{i,-1}(y)$ is 
constant on $(0,\delta)$.
\begin{comment}
We then observe that $\Delta_{t_1}f_{i,-1}(y)=c$ means $f_{i,-1}(y+t_1)-f_{i,-1}(y)=c$.
Let $y_1<\delta-t_1$, let $y_2=y_1+t_1$, and let $s>0$ be quite small compared to $t_1$ or $\delta$. 
We obtain 
$f_{i,-1}(y_2)=f_{i,-1}(y_1) +c$ and $f_{i,-1}(y_2+s)=f_{i,-1}(y_1+s) +c$
for all $y_1 \in (0,\delta-t_1)$ and $s>0$.
We observe $f_{i,-1}(y_2+s)-f_{i,-1}(y_2)=f_{i,-1}(y_1+s)-f_{i,-1}(y_1)$ for such $y_1$, $y_2$, and $s$.
Letting $s$ and $y_1$ vary, we conclude that for a dense subset of $(0,\delta)$, the set $(y, f_{i,-1}(y))$ is colinear with slope $\frac{c}{t_1}$.
By the assumption that $f_{i,-1}$ is continuous on $(0,\delta)$ and $\calR$ is o-minimal, and setting $b:= \lim_{y \to 0^+} f_{i,-1}(y)$, we can conclude that $f_{i,-1}(y)=b+\frac{c}{t_1}y$ on $(0,\delta)$.
\end{comment}
We observe that by Lemma \ref{delta}, the vector-valued function $ \hat{F}^{*}_1(y)=(f^*_1,\ldots,f^*_2)$ where $ f^{*}_1(y) = f_{i, -1}(y) \iff i=1 \text{ or }i \in J$, and $ f^{*}_i(y) = 0$ if $i \in [n] \setminus J$, and $ f^{*}_{n+1}(y) = h(y,f^{*}_2(y), \ldots ,  f^{*}_n(y))$,
coincides with the graph of an endomorphism in each pair of coordinates $(f^{*}_1, f^{*}_i)$ where $i \in J$ or $i = n+1$.
We iterate this process with $f_{i,-k}$ instead of $f_{i,-1}$ for each $k\in J$ to obtain a function $ \hat{F}^{*}_k(y)$ analogous to the function $\hat{F}_k$ described above.

As above, define $\tilde{h}^*_k( \hat{z}_k) = h(z_1, \ldots ,z_{k-1},0,z_{k+1},\ldots ,z_{|J|})$ where $z_i := f_{i,-k} (y)$ for the $i^{th}$ element of $J \setminus \{k\}$.
By UEP, there are finitely many possible values the coordinate functions $\tilde{h}^*_1, \ldots ,\tilde{h}^*_{|J|}$ can have on $(0,\delta)$.
By Lemma \ref{endoFact}, an EIC transformation is uniquely determined by how it acts on each coordinate, so we can recover a finite set of EIC transformations $h^*$ from the coordinate functions $\tilde{h}^*_1, \ldots ,\tilde{h}^*_{|J|}$.
We observe that we can write $\tilde{F}_k$ as a direct sum of $\hat{F}^*_k$ and $\pi^n_J(\tilde{F}_K)$.

By induction, we conclude that $\pi^n_J(\tilde{F}_K)$ only coincides with finitely many definable EIC transformations  as we range over the parameters.
As observed above, the same is true for $\hat{F}^*_k$.
Since any definable EIC transformation that coincides with $\tilde{F}_k$ on a neighborhood of the origin projects down to definable EIC transformations on $\hat{F}^*_k$ and $\pi^n_J(\tilde{F}_K)$, and any definable EIC transformation that coincides with $\tilde{F}_k$ can be uniquely recovered from those they induce on the direct summands, we conclude that there also are only finitely many possible EIC transformations that $h$ can be. 
This finishes the induction step.

\begin{comment}
For each iteration, we apply the induction hypothesis to the projection of $\hat{F}(y)$ onto the coordinates which are not identically zero on some interval with left endpoint $0$.
Define $h_k(\vec{x})$ to be the $n^{th}$ coordinate of $\hat{F}(y)$ for the $k^{th}$ iteration of this process.
We know that each $h_k(\vec{x})$ is a definable EIC transformation that is obtained from $h(\vec{x})$ by setting $x_i = 0$ for all $i$ in some subset $J_k \subseteq [n]$, with $|J_k| \geq 2$.
By induction, there are only finitely many definable EIC transformations that each $h_k$ can be.
Using a generalized version of Corollary \ref{endoFact}, we know we can recover $h$ uniquely from the $h_k$'s,
hence there also are only finitely many possible EIC transformations that $h$ can be, finishing the induction step.
\end{comment}
\qedhere
\end{proof}
\end{lem}

\begin{cor}\label{fincor}
Suppose that $F(\vec{x},\vec{y}): U \times A \subseteq M^{n+k} \to M^m$ is $\calL (\emptyset)$-definable, and for every $\vec{a} \in A \subseteq M^k$ (the parameter space) $F(\vec{x}, \vec{a})$ is a continuous function with domain $U \subseteq M^n$.
Then 
\begin{gather*}
\{ \hat{F}(\vec{x}) = [ F(\vec{x} +\vec{b}, \vec{a})-F(\vec{b},\vec{a}) ]|_{B^n_{\epsilon}(0)} : \vec{a} \in A, \vec{b} \in U, \epsilon >0, \gr( \hat{F}) \subseteq H \subseteq M^m 
 \text{a definable hyperplane}\}
\end{gather*}
collapses to a finite definable family of germs of functions (as in Definition \ref{germs}) through the origin.
\begin{proof}
We proceed by induction on $n = |\vec{x}|$.
For the base case, suppose that $\vec{x} = x$ is a single variable.
We consider the following family of functions with restricted domain:
$$\calF :=\{\hat{F}(x +b, \vec{a}):=F(x +b, \vec{a})-F(b,\vec{a})|_{B_{\epsilon}(0)} : b \in U, \vec{a} \in A, \epsilon >0, \gr(\hat{F}) \subseteq H \text{ a def. hyperplane}\}$$
and consider the set of $\vec{a} \in A$ and $b \in M$ such that $\exists \delta >0$ for which $\forall x,y \in B_{\delta}(0) (\hat{F}(x+y+b,\vec{a}) = \hat{F}(x+b, \vec{a}) + \hat{F}(y+b, \vec{a}))$. 
This carves out a $\emptyset$-definable subfamily of $\calF$ that collapses (as in Definition \ref{germs}) to a collection of germs of endomorphisms in each coordinate. 
We apply Lemma \ref{demarcation} to conclude that for all $b \in M$ and $\vec{a} \in A$ there are only finitely many definable hyperplanes which coincide with $\hat{F}(x+b, \vec{a})$ on some neighborhood.

Now let $n>1$ and assume the claim holds for $n-1$.
Write $\vec{x} = (x_1, \ldots , x_n)$, and let 
$\mathcal{H} = \{ H_{\vec{b},\vec{a}} \subseteq M^m : \vec{a} \in A, \vec{b} \in U \}$ enumerate the definable hyperplanes $H_{\vec{b},\vec{a}}$ that coincide with $\hat{F}(\vec{x} +\vec{b}, \vec{a})$ for some parameters $\vec{a} \in A$ and $\vec{b} \in U$. 
Since Lemma \ref{demarcation} holds for every arity of parameter tuple, we simply ``move'' the last variable $x_n$ from the domain of the function to the parameter space. 
By this, we mean that we can think of the $n$-dimensional hypersurface $\hat{F}(\vec{x} + \vec{b}, \vec{a})$ as a definable family of $n-1$-dimensional hypersurfaces, i.e. $\{ \hat{F}( x_1+b_1, \ldots ,x_{n-1}+b_{n-1}, c, \vec{a}):M^{n-1}\to M^m :\vec{b} \in \pi^{[n-1]}(U), c \in \pi^{n}(U), \vec{a} \in A \}$, in which $(c, \vec{a})$ is now the parameter tuple ranging over $\pi^n(U) \times A$, with domain $\pi_{n-1}^{n} (U)$.
We will write $\hat{x}_n = (x_1, \ldots ,x_{n-1})$.
We now use the induction hypothesis to conclude there are only finitely many definable hyperplanes, say $H_1, \ldots ,H_{\ell} \subseteq M^{m} $, for which there is some $\epsilon$-neighborhood of zero and some $\hat{F}(\hat{x}_n +\hat{b}_n, c, \vec{a})$ in this family of $n-1$-dimensional hypersurfaces that coincides with some $H_i$ on said $\epsilon$-neighborhood.

Finally, we appeal to the fact that the set of all $m$-dimensional hyperplanes that contain some level set of a hypersurface (intersected with a neighborhood of zero) is by definition a superset of the $m$-dimensional hyperplanes that contain the entire surface (intersected with the same neighborhood of zero).
This can easily be seen by writing out the set $\mathcal{H}$ for the family $F$, and the set $\mathcal{H}$ defined analogously for the family $F(\hat{x}_n,z, \vec{y})$ (in which the $n^{th}$ coordinate is regarded as a parameter) as follows:
\begin{gather*}
\mathcal{H}_F = \{ \hat{F}(\vec{x} +\vec{b}, \vec{a})|_{B^n_{\epsilon}(0)} : \vec{a} \in A, \vec{b} \in U, \epsilon >0, \gr( \hat{F}) \subseteq H \subseteq M^m \text{ a definable hyperplane} \} \subseteq \\
\{ \hat{F}(\hat{x}_n +\hat{b}_n, c,\vec{a})|_{B^{n-1}_{\epsilon}(0)} : c \in \pi^n(U), 
\vec{a} \in A, \vec{b} \in U, \epsilon >0, \gr( \hat{F}) \subseteq H \subseteq M^m \text{ a def. hyperplane} \}\\
 = \hat{\mathcal{H}}_F.
\end{gather*}
Viewing $\mathcal{H}_F$ as a collection of germs, we see that it is a finite collection because every hyperplane $H \in \mathcal{H}_F$ that coincides with $\hat{F}(\vec{x},\vec{a})$ on a neighborhood $U'$ also coincides with the level sub-hypersurface $\hat{F}( \hat{x}_n, c, \vec{a})$ for any $c \in \pi^n(U')$. 
This concludes the induction argument.
\end{proof}
\end{cor}

By restricting our attention to o-minimal theories $T$ which do not have UEP, we can ensure by Lemma \ref{demarcation} that we can axiomatize the following property: given a generic input for a definable function $f:M^n \to M^m$ that is not a definable EIC transformation, we can obtain points in the image of $f$ both inside and outside of the predicate group.

\begin{defn} \label{cfdef}
Let $(\calM, \calG) \models T_{\calG}$ and $n \in \N_{>0}$. 
Let $\vec{f} = (f_1, \ldots f_k):U \subseteq M^n \to M^k$ be an $n$-ary vector function that is $\calL(\emptyset)$-definable.
For $I \subseteq [k]$, let $M^I$ denote $\prod_{i = 1}^k P_i$ where $P_i = M$ if $i \in I$ and $P_i = \{0\}$ if $i \not \in I$.
Given $\vec{d} \in U$ and $I \subseteq [k]$, let 
$$\tilde{f}^{d,I}_i(\vec{x}) = \begin{cases} f_i(\vec{x}) - f_i(\vec{d}), & i \in I \\ 0,& i \not \in I \end{cases}.$$

We let $\mathcal{Q}$ denote the (non-definable) set of definable hyperplanes $H \subseteq M^{n+k}$ that have a matrix representation with exclusively $\Q$-affine entries.
Let $\tilde{f}^{d,I} = (\tilde{f}_1, \ldots , \tilde{f}_k)$.
We define the following:
\begin{gather*}
C_f = \Big \{ \vec{d}=(d_1, \ldots ,d_n) \in U : \text{ for some } \emptyset \neq I \subseteq [k] \text{ there is } H \subseteq M^{n+k} \text{ with } H \in \mathcal{Q}\\
\text{and } \exists \epsilon > 0 \forall x \in B_{\epsilon}(\vec{d}) \left( (\vec{x},\tilde{f}^{d,I}(\vec{x})) \in H \right) \Big \}.
\end{gather*}
Furthermore, let $\tilde{C}_f$ denote $\interior (U \setminus C_f)$, the topological interior of the complement of $C_f$ in $U$.
\end{defn}

We see that $C_f$ is open because by definition the property holds for an element only if it holds on a neighborhood of that element, and it is a subset of an open set.
Intuitively, we think of $C_f$ as being the open subset of the domain of $\vec{f}$ on which $\vec{f}$ or some projection of $\vec{f}$ onto a subspace of $M^n$ is locally $\Q$-affine.
We think of $\tilde{C}_f$ as the obstacle-free part of the domain of a definable function.
The following remark is really a corollary to Lemma \ref{demarcation}.
\begin{rmk} \label{cfrmk}
For all $\vec{f}$ as described in the above definition, the set $C_f$ is $\calL$-definable.
\begin{comment}
\begin{proof}
This follows immediately from the assumption that $T$ does not have UEP in conjunction with Lemma \ref{matrix} and Lemma \ref{demarcation}.
These ensure there are only finitely many definable hyperplanes in $\mathcal{Q}$ (as defined above) with which the graph of some sub-tuple of $\vec{f}$ coincides; enumerate these hyperplanes as $\{H_1, \ldots ,H_t\} \subseteq \mathcal{Q}$.
We consider $S = S_1 \cup \cdots \cup S_t$
where $S_i$ is the subset of the domain of $\vec{f}$ on which the graph of some sub-tuple of $\vec{f}$ agrees with $H_i$.
Realizing that the set $S$ is exactly the set $C_f$, we conclude we can write it as a finite union of $\calL( \emptyset )$-definable sets.
\end{proof}
\end{comment}
\end{rmk}

By the $\calL$-definability of $C_f$ and the fact that $T$ is a complete theory, we conclude that for $\calL( \emptyset)$-definable sets and functions, the set $C_f$ should have the same first-order properties with respect $\calL$ in every model of $T_{\calG}$.
Below, for $\vec{x}=(x_1, \ldots ,x_n) \in M^n$ we will use $\vec{x}_{\iota}$ to denote $(x_{i_1}, \ldots ,x_{i_k})$ where $\iota = (i_1, \ldots ,i_k) \in [n]^k$. 
Let $[n]_k$ denote the set of all $k$-element subsets of $n$.

We now show in Lemma \ref{suitable} that we may apply Proposition 1.12 from \cite{E18} to the theory $T_{\calG}$ when $T$ does not have UEP.
Below $T_0$ is a reduct of $T$ with sub-language $\calL_0 \subseteq \calL$ such that the algebraic closure operation (call it $\acl_0$) is a pregeometry, and $\Ofork$ is the independence relation associated to $\acl_0$.
D'Elb\'{e}e uses the following definition in his companionability criterion, in which $T_S$ is the theory of pairs $(\calM, \calM_0)$ where $\calM \models T$ and $\calM_0 \models T_0$.
\begin{defn*}[\cite{E18}, 1.10]
We say that a triple $(T,T_0,\calL_0)$ is \emph{suitable} if it satisfies the following:
\begin{itemize}
\item[($H_1$)] $T$ is model complete;
\item[($H_2$)] $T_0$ is model complete, and for all infinite $A$, $\acl_0(A) \models T_0$;
\item[($H_3^{+}$)] $\acl_0$ defines a modular pregeometry;
\item[($H_4$)] for all $\calL$-formula $\vPhi(x,y)$ there is some $\calL$-formula $\theta_{\vPhi}(y)$ such that for $b \in \calM \models T$, 
\begin{align*} 
\calM \models \theta_{\vPhi}(b) \iff & \text{ there exists } \calN \succcurlyeq \calM \text{ and }a \in \calN \text{such that } \\  & \vPhi(a,b) \text{ and }a\text{ is } \Ofork \text{-independent over } \calM .
\end{align*}
\end{itemize}
\end{defn*}

Below, the theory $TS$ mentioned in the criterion D'Elb\'{e}e gives for companionability is axiomatized in \cite{E18} as follows. 
For $x = x^0x^1$, and for each $\calL$-formula $\vPhi(x,y)$ and each $(\tau_i(t,x,y))_{i<k}$ a finite set of $\calL_0$-formulas that are algebraic in $t$ and strict in $x^1$ (i.e. if we vary $x^1$ then $t$ varies as well), the sentence
\[ \forall y \Bigg ( \theta_{\vPhi}(y) \rightarrow \Big ( \exists x \vPhi (x,y) \land x^0 \subseteq S \land \bigwedge_{i<k} \forall t(\tau_i(t,x^0,x^1) \rightarrow t \not \in S ) \Big ) \Bigg) \]
is in $TS$, and the set of all such sentences along with those in $T$ axiomatize $TS$.
In our setting, we will call these axioms the ``companion axioms.''

\begin{prop*}[\cite{E18}, 1.12]
Let $(T,T_0,\calL_0)$ be a suitable triple. Then $TS$ exists and is the model-companion for the theory $T_S$.
\end{prop*}

We now define a triple $(T,T_0,(S_1,S_2;0,+,(q\cdot(-))_{q\in \Q}, (c_i)_{i \in \N}))$ that will be our candidate for such a suitable triple, and whose theory is bi-interpretable with the theory $T_{\cal{G}}$.
Moreover, in the case that $T$ has UEP, the model companion associated with this triple (which we will shortly conclude is given to us by Proposition 1.12 above) will also have a model companion axiomatized analogously via the evident bi-interpretation.
For the theory $T_{\calG}$, recall that ${\bf D}$ is the open, connected, and unary $\calL(\emptyset)$-definable superset of $\calG$ in which $\calG$ is dense and codense.
Let $T_0$ have models given by $(S_1,S_2; 0,+,(q\cdot(-))_{q\in \Q}))$, where $S_1$ and $S_2$ are the two sorts of $T_0$, such that $S_1=\mathbf{D}$ and $S_2=\mathbf{D}^c$, and $T_0$ dictates that $(S_1, 0,+,(q\cdot(-))_{q\in \Q})) \models \operatorname{ODAG}$ and the structure on $(S_2, (c_i)_{i \in \N})$ is simply that of an infinite set in which all $c_i$'s are distinct.
 
\begin{lem}\label{suitable}
If $T$ does not have UEP, then the triple $(T,T_0,(S_1,S_2;0,+,(q\cdot(-))_{q\in \Q}, (c_i)_{i \in \N}))$ is suitable.
\begin{proof}
We already require ($H_1$) of $T$.
In the language $(0,+,q\cdot (-))_{q \in \Q})$ the theory ODAG is model complete and has quantifier elimination, and similarly the theory of infinite sets has quantifier elimination in the language $(c_i)_{i\in \N}$ of infinitely many constants.
Hence it is evident that $T_0$ is model complete in the language $\calL_0:=(S_1,S_2;(q\cdot(-))_{q\in \Q}, (c_i)_{i \in \N})$, as ($H_2$) requires.
Given a set $A$ in a model of $T$, the $\calL_0$-algebraic closure $\acl_0(A)$ is the disjoint union of the $\Q$-linear span of $A$ (which is a model of ODAG) with the set $\{ c_i:i \in \N\}$, so condition ($H_2$) is satisfied.
Moreover, since $\acl_0$ is the disjoint union of the $\Q$-linear closure of a set with $\{ c_i:i \in \N\}$, the pregeometry it defines is modular, and hence ($H_3^+$) is satisfied.

To see that condition ($H_4$) holds, suppose that $\vPhi (\vec{x},\vec{z})$ is a $\calL$-formula with $|\vec{x}|=d$ and $|\vec{z}| = m$.
Let $\mathcal{Z} := \{ \vec{z}: \exists^{\infty} \vec{x}\vPhi(\vec{x},\vec{z})\}$ and let $\mathcal{Y} := \{(\vec{x},\vec{z}): \vec{z} \in \mathcal{Z} \land \vPhi(\vec{x},\vec{z})\}$.
There are two cases; the first is that $\mathcal{Y} \cap (\mathbf{D}^c)^n$ is nonempty.
In this case, the pregeomtery induced by $\acl_0$ on $\mathcal{Y} \cap (\mathbf{D}^c)^n$ is trivial, so by compactness we can always find some $\calN \succcurlyeq \calM$ such that $\mathcal{Y} \cap (\mathbf{D}^c)^n \setminus \{c_i: i\in \N\}$ is nonempty and any element of this set is $ \Ofork$-independent over $\calM$.
The second case is that $\mathcal{Y} \cap (\mathbf{D}^c)^n$ is empty, in which case we proceed as follows.
Note that if $\pi_i(\mathcal{Y})$ is contained in $\mathbf{D}^c$ for any coordinate $i \in [d]$, any choice for $x_i$ from $\mathbf{D}^c \setminus \{c_i:i \in \N\}^{\calN}$ will be $\Ofork$-independent over $\calM$ from the rest of the tuple so long as $\calN \succcurlyeq \calM$ is at least countably saturated.

\begin{comment}
We partition $\mathcal{Y}$ into definable sets $D_1, \ldots ,D_d$ and $E$ such that the following holds.
For each $n \in [d]$ and all $(\vec{x},\vec{z}) \in D_n$, $n$ is the unique element of $[d]$ for which there is a sub-tuple $\vec{x}_{[n]}$ of $\vec{x}$ and there are definable functions $\vec{f}\,^1, \ldots ,\vec{f}\,^{k_n}$ such that $\vec{f}\,^i(\vec{x}_{[n]}, \vec{z}) = (f_{n+1}(\vec{x}_{[n]}, \vec{z}), \ldots ,f_{d}(\vec{x}_{[n]},\vec{z}))$ for each $i \in [k_n]$ and a definable partition of $D_n$ into $D^1_n, \ldots ,D^{k_n}_n$  such that 
$$T \models \forall (\vec{x},\vec{z}) \in D^i_n \Bigg ( \Big ( \vPhi(\vec{x},\vec{z}) \leftrightarrow \vPhi(\vec{x}_{[n]},\vec{f}\,^i(\vec{x}_{[n]},\vec{z}),\vec{z})\Big ) \land \vec{x}_{[n]} \not \in C_{\vec{f}^i}\Bigg).$$
Recall the $\calL(\emptyset)$-definability of the set $C_f$ as defined in Definition \ref{cfdef}.
Additionally, for each $n \in [d]$ and each $i \in [k_n]$, we require that for each $\vec{z} \in D_n$ the set $\{\vec{x}: (\vec{x},\vec{z}) \in D^i_n\}$ either is empty or has dimension exactly $n$, and we also require that for each such $\vec{z}$ the projection of $D_n$ onto the sub-tuple $\vec{x}_{[n]}$ is open.
We may take each $D_n$ to be the maximal subset of $\mathcal{Y}$ on which the above holds, and we also may assume each $k_n$ could not be made smaller without violating one of the above properties.
Let $E$ contain all the tuples $(\vec{x},\vec{z})$ in  $\mathcal{Y} \setminus D_1 \cup \ldots \cup D_d$.
\end{comment}

As a definable set in an o-minimal structure, $\mathcal{Y}$ has a cell decomposition of the following form:
$$\{D^j_i: i \in [d], j\in [k_i], k_1, \ldots k_d \in \N \}$$ 
where for each $1\leq i \leq d$ and for each $1\leq j \leq k_i$, the set $\{\vec{x}: \exists \vec{z}((\vec{x},\vec{z})\in D^j_i)\}$ has natural dimension $i$, and there are precisely $k_i$ such cells.
Given $n \in [d]$, for each $j \in [k_n]$ and for all $(\vec{x},\vec{z}) \in D^j_n$, let $\vec{x}_{j,[n]}$ denote the subtuple of $\vec{x}$ such that for all $(\vec{x},\vec{z}) \in D^j_n$, the projection on $D^j_n$ onto the $n$ coordinates of the subtuple $\vec{x}_{j,[n]} \subseteq \vec{x}$ is an open subset of $n$-space.
Permuting the order of the coordinates of $\vec{x}$ in the definition of $D^j_n$ if necessary, we may assume without loss of generality that $\vec{x}_{j,[n]} = (x_1,\ldots ,x_n)=\vec{x}_{[n]}$.
For each $n \in [d]$ and $j \in [k_n]$, the above cell decomposition give us vector functions of the following form:$$\vec{f}^{j}(\vec{x}^j_{[n]},\vec{z})=(f^j_{n+1}(\vec{x}^j_{[n]},\vec{z}), \ldots ,f^j_{d}(\vec{x}^j_{[n]},\vec{z}))$$
such that the following holds:
$$\forall (\vec{x},\vec{z}) \in D^j_n \Big ( \vPhi(\vec{x},\vec{z}) \leftrightarrow \vPhi(\vec{x}_{[n]},\vec{f}\,^j(\vec{x}_{[n]},\vec{z}),\vec{z})\Big ).$$

Recall the $\calL(\emptyset)$-definability of the set $C_f$ as defined in Definition \ref{cfdef}.
Let $D_n$ be the union of the following sets:
$$\{(\vec{x},\vec{z}): (\vec{x},\vec{z})\in D^j_n \land \vec{x}_{[n]} \in \tilde{C}_{\vec{f}^j} \}$$
for all $j \in [k_n]$.
We think of $D_n$ as being the intersection of the complement of $\{C_{\vec{f}^j}:j\in [k_n]\}$ with the subset of the projection of $\vPhi$ onto $\vec{x}$ that has natural dimension $n$.
Let $E$ contain all the tuples $(\vec{x},\vec{z})$ in $\mathcal{Y} \setminus D_1 \cup \ldots \cup D_d$.
Consider the following formula:
\[ \theta_{\vPhi} (\vec{z}) := \bigvee_{n=1}^d \exists \vec{x} \vPhi(\vec{x},\vec{z}) \land (\vec{x},\vec{z}) \in D_n. \]
We now show that $\theta_{\vPhi}(\vec{z})$ holds precisely if there is $\calN \succcurlyeq \calM $ and $\vec{a}$ in $\calN$ such that $\calN \models \vPhi(\vec{a},\vec{b})$ and $a$ is $\Ofork$-independent over $\calM$. 

First, we suppose that $\calM \models T$, $\vec{b} \in M^{|\vec{z}|}$ and $\calM \models \theta_{\vPhi}(\vec{b})$.
Let $\calN \succcurlyeq \calM$ be a $\calM^{+}$-saturated elementary superstructure (in particular, also $\max \{ | \calL |^{+}, \aleph_1 \} $-saturated).
Let $n \leq d$ be such that $(\vec{x}, \vec{b}) \in D_n$ for some $\vec{x} \in M^n$, and suppose $(\vec{x},\vec{z}) \in D^i_n$ in the cell decomposition of $\mathcal{Y}$.
Since the projection onto the first $n$ coordinates of $D_n$ is open (by definition) we know that the partial type over $M$ given by
$$\Theta_{\vec{b}} = \Bigg \{ \vPhi(\vec{x}_{[n]}, \vec{f}\,^i(\vec{x}_{[n]},\vec{b}\,),\vec{b}\,) \land \Big( \vec{x}_{[n]} \not  \in C_{\vec{f}^i} \Big) \Bigg \} \bigcup \Bigg \{ \vec{q} \cdot (\vec{x}_{[n]},\vec{f}\,^i(\vec{x}_{[n]},\vec{b}\,)) \neq m : \vec{q} \in \Q^{d}, m \in \calM \Bigg \} $$
is consistent with the theory of $\calM$ because it is finitely satisfiable, by the definition of $C_{\vec{f}}$.

Hence $\Theta_{\vec{b}}$ has a realization $\vec{a} \in N^d$.
It is clear that if $\vec{a}$ satisfies $\vPhi(\vec{x},\vec{b})$ and the sub-tuple $\vec{a}_{[n]}$ satisfies the partial type $\Theta_{\vec{b}}$, 
then it does not lie in any set that is $\acl_0$-dependent over $\calM$, since those are precisely the $\Q$-affine hyperplanes defined over $M$.
So the ($\implies$) part of condition ($H_4$) holds.

Now suppose that there is $\calN \succcurlyeq \calM $ and $\vec{a} \in \calN$ such that $\calN \models \vPhi(\vec{a},\vec{b}\,)$ and $\vec{a}$ is $\Ofork$-independent over $\calM$. 
Note that if $(\vec{a},\vec{b}\,) \not \in \mathcal{Y}$, then $\vec{a}$ would have to be in the $\dcl$ of $\vec{b}$
in which case $\vec{a} \in \calM$, which would contradict the $\Ofork$-independence of $\vec{a}$ over $\calM$.
If $(\vec{a},\vec{b}\,) \in E$, due to o-minimal cell decomposition either one of two things must hold.
The first possibility is that for some $n \in [d]$ and $i \in [k_n]$ we must have $\vec{a}_{[n]} \in C_{\vec{f}^i}$.
Yet by definition of $C_{\vec{f}^i}$ it would then follow that $\vec{a}$ is not $\Ofork$-independent over $\calM$.
The other possibility is that $\vec{a} \notin C_{\vec{f}^i}$ for some $n \in [d]$ and $i \in [k_n]$, yet the projection of the following onto the first $n$ coordinates is not open for each $\epsilon>0$:
$$\{ \vec{x} \in B_{\epsilon}(\vec{a}): \vPhi(\vec{x}_{[n]},\vec{f} (\vec{x}_{[n]},\vec{b}),\vec{b}) \land \vec{x} \not \in C_{\vec{f}})\}.$$
If this is the case then for all $\epsilon>0$, the set of $\vec{x}=(\vec{x}_{[n]}, \vec{f}^i (\vec{x}_{[n]},\vec{b})) \in B_{\epsilon}(\vec{a})$ outside of a finite set of hyperplanes dictated by $C_{\vec{f}^i}$ is itself contained in a finite set of hypersurfaces, due to cell decomposition. 
However, the open ball $B_{\epsilon}(\vec{a})$ must intersect infinitely many hyperplanes and hypersufraces, so it is impossible to have $\vec{a} \notin C_{\vec{f}^i}$ and not also have $B_{\epsilon}(\vec{a}) \not \subseteq C^c_{\vec{f}^i}$.

Hence $\vec{a}$ witnesses also that $\vec{b}$ satisfies $\theta_{\varphi}$, as desired.
Thus $\calM \models \theta_{\vPhi}(\vec{b})$, and condition ($H_4$) holds, as desired.
\end{proof}
\end{lem}

As a corollary to Lemma \ref{suitable}, and Proposition 1.12 in \cite{E18}, we are ready to conclude Theorem \ref{EWconjecture} holds.
We define a slight variant of the formula $\theta_{\vPhi}$ from before:
\[\tilde{\theta}_{\vPhi} (\vec{z}) := \bigvee_{n=1}^d \bigvee_{i=1}^{k_n} \exists \vec{x} \vPhi(\vec{x},\vec{z}) \land \vec{x} \in {\bf D}
\land (\vec{x},\vec{z}) \in D^i_n. \]

\begin{thm}[Theorem \ref{EWconjecture}]\label{gen dichotomy}
For $T$ an o-minimal extension of the theory of ordered abelian groups with language $\calL$, let $T_{\calG}$ be the $\calL \cup \{ \calG \}$-theory extending $T$ that states also that $\calG$ is a divisible subgroup both dense and codense in ${\bf D}$. 
If the only families of additive endomorphisms whose graphs are uniformly definable in a neighborhood of $0$ have finitely many germs at $0$, then $T_{\calG}$ has a model companion $T^*_{\calG}$ given as follows.
For $\vec{x} = \vec{x}_0\vec{x}_1$, and for each $\calL$-formula $\vPhi(\vec{x},\vec{y})$ and each $(\tau_i(\vec{t},\vec{x},\vec{y}))_{i<k}$ a finite set of $\calL_0$-formulas that are algebraic in $\vec{t}$ and strict in $\vec{x}_1$ (i.e. if we vary $\vec{x}_1$ then $t$ varies as well), the sentence
\[ \forall \vec{y} \Bigg (\tilde{\theta}_{\vPhi}(\vec{y}) \rightarrow \Big ( \exists \vec{x} \vPhi (\vec{x},\vec{y}) \land \vec{x}_0 \subseteq \calG \land \bigwedge_{i<k} \forall \vec{t}\Big( \tau_i ( \vec{t},\vec{x}_0,\vec{x}_1) \rightarrow t \subseteq \calG^c \Big) \Big) \Bigg) \]
is in $T^*_{\calG}$, and $T_{\calG} \subseteq T^*_{\calG}$.
Otherwise, $T_{\calG}$ does not have a model companion.

\begin{proof}
Under the given hypotheses, the axiomatization for $T^*_{\calG}$ follows from Proposition 1.12 in \cite{E18}. 
We can apply Proposition 1.12 since Lemma \ref{suitable} shows the hypotheses of the proposition are satisfied if $T$ does not have UEP.
The proof of the ``otherwise'' statement is given by Theorem \ref{negprop}.
\end{proof}
\end{thm}

We will see in the examples section below that this characterization translates to an even cleaner dichotomy when applying the characterization to an o-minimal expansion $\calR$ of a real closed field with an added predicate for a multiplicative subgroup dense in $\calR_{>0}$.

\section{Examples}\label{ex}

We consider the connections that the results in this paper have to the framework of \cite{BGHK18}.
Let $\R_K$ denote $(\R, <, +, (x \mapsto kx )_{k \in K} )$, where $K \subseteq \R$ is a subfield.  
We observe that the structure  $(\calR, \calQ) = (\R, <, +, (x \mapsto kx )_{k \in K}, \Q )$ is a model of $T_{\calG}$. The theory of this structure is not model complete, though $T_{\calG}$ in this case does have a model companion.
In the language $\calL_U = (0,1, <, +, ( k(x))_{k \in K},U(x))$, where $k(x)$ is the symbol for scalar multiplication by $k \in K$ and $U$ is the predicate that picks out $\calQ$, we know that $T^d_K := Th( \R_K, \Q)$ is what the authors of \cite{BGHK18} call an ``ML theory.''
In particular, this implies that it is near-model complete.  Yet it is not model complete, and what fails is linear disjointness.

\begin{ex}\label{VSex}
Consider the field $K = \R^{alg} (e)$.  Let $\calR_1 = \R^{alg}( e, \zeta, \eta)$ where $\zeta$ is a transcendental number over $K$, and $\eta$ is a transcendental number over $K( \zeta)$.  Let $\calR_0 = \R^{alg}(e, \zeta)$.  Let $\calQ_0 = \Q $, and let $\calQ_1 = \Q ( -e \eta + \zeta , \eta)$.
Then it is easy to check $(\calR_1, \calQ_1)$ and $(\calR_0, \calQ_0)$ are models of $T^d_K$, and by construction $(\calR_0, \calQ_0) \subseteq ( \calR_1, \calQ_1)$.  
However it is not the case that $(\calR_0, \calQ_0) \preccurlyeq (\calR_1, \calQ_1)$ since $(\calR_1, \calQ_1) \models \exists q_1 \exists q_2 U(q_1) \land U( q_2) \land ( e(q_1) + q_2 = \zeta)$, but $( \calR_0, \calQ_0) \models \forall q_1, q_2 (U(q_1) \land U(q_2) \rightarrow e(q_1) + q_2 \neq \zeta$.
\end{ex}

\subsection{The case that $+$ is addition in a field.}

We note that by Theorem \ref{gen dichotomy}, when $T$ extends the theory of real closed fields and $\calG$ picks out an additive subgroup, the theory $T_{\calG}$ does not have a model companion.
However, this is also observed by d'Elb\'ee in \cite{E18}; see Proposition 5.32 and Remark 5.33 in \cite{E18}.

\subsection{The case that $+$ is multiplication in a field.}
We can apply the companionability characterization to pairs in which the underlying o-minimal structure $\calM$ expands a real closed field, and the group $\calG$ is a multiplicative subgroup of $\calM_{>0}$.
Interpreting $0, +$ in $\calL$ from section \ref{(R,G)} as $1^{\calM}$ and $\cdot^{\calM}$, we can clearly establish the same companionability dichotomy for o-minimal extensions of RCF. We make this more concrete in the examples below.

\begin{ex}\label{mult non-ex}
An example of a non-companionable structure is $(\R_{\eX}, G)$ where $G < \R_{>0}$ is a dense, divisible multiplicative subgroup. 
This follows since $\ln (x)$ is definable in $\R_{\eX}$, so $y^x = e^{x\ln (y)}$ is definable. 
Hence the family of multiplicative endomorphisms $\{ x^r : 0<r<1 \}$ is uniformly definable in $\R_{\eX}$, and witnesses that $Th(\R_{\eX})$ has UEP.
By Theorem \ref{negprop}, this implies no model companion exists.
\end{ex}

Recall that a partial exponential function on $\calM$ is a definable differentiable function 
$f: (-\epsilon, \epsilon)\to M_{>0}$ for some $\epsilon>0$ in $M$, such that $f’=f$. 

\begin{thm}[\ref{mult}]\label{mult dichotomy}
Let $T$ be an o-minimal theory extending the theory of real closed fields (RCF), and let $T^{\times}_{\calG}$ be the extension of $T$ to the language $\calL_{\calG}$ in which $\calG$ is a predicate for a divisible multiplicative subgroup of the positive elements that is dense and codense as a subset of an open unary set definable in $\calL$.
Then $T^{\times}_{\calG}$ has a model companion if and only if for every model $\calM \models T$ there is no partial exponential function that is $\calL$-definable in $\calM$.
\begin{proof}
$(\implies)$ If a partial exponential function $f$ is definable on some interval $I \subseteq M$, then we can apply Theorem \ref{negprop} to conclude that $T^{\times}_{\calG}$ does not have a model companion, as illustrated in Example \ref{mult non-ex}.

$(\impliedby)$ We suppose that $\calM$ defines an infinite family of local endomorphisms of the multiplicative group $M_{>0}$, and we show this implies that a function defined on some interval $I$ satisfies the differential equation that characterizes a partial exponential function.
We let $\tilde{F}(x, y)$ be any definable partial function in $\calM$ such that its domain is the box $I \times \tilde{J} \subseteq M_{>0} \times M$, where we assume  $1 \in I$, and $\tilde{F}(x_1,y)\tilde{F}(x_2,y) = \tilde{F}(x_1x_2,y)$ for all $x_1,x_2 \in I$ and $y \in \tilde{J}$.
Such a function $F$ can be defined in $\calM$ using the existence of an infinite definably family of functions, the o-minimality of the structure $\calM$, and the definable choice functions o-minimality gives us.

We write $\tilde{F}(x,y) = \tilde{f}_y(x)$ for a fixed $y \in \tilde{J}$. 
We now define $J = \{ \tilde{f}'_y(1) :y \in \tilde{J} \}$, where $\tilde{f}'_y(x) = \frac{d}{dx} \tilde{f}_y(x)$, and let $\eta(y) =\tilde{f}'_y(1)$, which we can without loss of generality take to be continuous and injective, restricting $\tilde{J}$ if needed.
One can show that by possibly inverting and shifting some of the functions $\tilde{f}_y$, we can make $0$ the left endpoint of $J$.
Let $F(x,y) : I \times J \to M$ be given by $F(x,y) = \tilde{F}(x, \eta^{-1}(y))$.
We now deduce that since $f_y(x_1)f_y(x_2) = f_y(x_1x_2)$ we also have $\frac{\partial}{\partial z} f_y(xz) = \frac{\partial}{\partial z} f_y(x)f_y(z) = f_y(x)f'_y(z)$ and $\frac{\partial}{\partial z} f_y(xz) = xf'_y(xz)$, 
so $f_y(x) = f'_y(xz) \cdot \frac{x}{f'_y(z)}$.
Letting $z=1$, we conclude that for each $y \in J$ we know $f_y(x) = f'_y(x) \cdot \frac{x}{f'_y(1)}$ characterizes the function $f_y$ subject to the constraints %$f_y(0) =0$ ??
$f_y(1)=1$ and $\forall x,z \in I( f_y(xz) = f_y(x)f_y(z))$.
We hence conclude for any $y_1,y_2 \in J$ that $f'_{y_1}(1) = f'_{y_2}(1) \iff f_{y_1} = f_{y_2}$.

We will use this to show that $F(x,y_1+y_2) = F(x,y_1) \cdot F(x,y_2)$.
As a product of definable endomorphisms, we observe that $F(x,y_1) \cdot F(x,y_2)$ is a definable endomorphism as well.
By our above remarks, $f'_{y_1+y_2}(1)$ uniquely determines the function $f_{y_1+y_2}$ subject to the constraint of being a definable endomorphism.
We know that $f'_{y_1+y_2}(1) = y_1+y_2$ by how we defined $F(x,y)$ and $J$ with the use of $\eta(y)$.
We also observe that $\frac{\partial}{\partial x} F(x,y_1)F(x,y_2) = F(x,y_1)F'(x,y_2) + F(x,y_2)F'(x,y_1)$.
Since $f_y(x)$ is an endomorphism on $I$, we know that $f_y(1)=1$ for all $y \in J$.
So $\frac{\partial}{\partial x} F(x,y_1)F(x,y_2)$ at $x=1$ is $F(1,y_1)F'(1,y_2) + F(1,y_2)F'(1,y_1) = y_2 + y_1$,
hence $f_{ y_1} (x)f_{ y_2 } (x)$ has $y_1+y_2 = f'_{ y_1+y_2 }(1)$ as its derivative at $x=1$, which makes them the same function on $I$.

Finally, we show that $\frac{\partial}{\partial y} F(x,y) = C(x) F(x,y)$ where $C(x)$ is purely a function in $x$ to conclude that $F(e,y)$ is equivalent to an exponential function (or shift thereof) on some subinterval of $J$.
Observe that 
$$\lim_{h \to 0} \frac{F(x,y+h) - F(x,y)}{h} = \lim_{h \to 0} \frac{F(x,y)F(x,h) - F(x,y)}{h} = F(x,y) \lim_{h \to 0} \frac{F(x,h) - 1}{h}$$
so letting $C(x) = \lim_{h \to 0} \frac{F(x,h) - 1}{h}$ we obtain the desired result.
Since the differential equation $\frac{\partial}{\partial y} F(e,y) = C(e) F(e,y)$ characterizes all functions of the form $c_1e^{c_2y} + c_0$ where $c_0,c_1,c_2 \in \R$, we conclude that the function $F(e,y)$ coincides with a partial exponential function on a subinterval of $J$.
\end{proof}
\end{thm}

When $T$ is an o-minimal extension of $RCF$ and $T^{\times}_{\calG}$ is the extension described in Theorem \ref{mult}, we use $T^*_{\calG}$ to denote the model companion of $T^{\times}_{\calG}$, if it exists.
Depending on the theory $T$, it will be clear from context which kind of model companion $T^*_{\calG}$ denotes for the remainder of this paper.

\begin{ex}
For the structure $\calR^* := (\R, 0,1, +, \cdot, (k)_{k \in K}, (x^k)_{k \in K})$, where $K \subseteq \R$ is any subfield, the pair $(\calR^*, \calG)$, with $\calG \subseteq \R_{>0}$ dense and codense in $\calR_{>0}$, has a model companion.
\begin{proof}
If $r^x$ is definable on some interval in $\calR$, then it is definable from the real field $\bar{\R}$ augmented with finitely many functions of the form $x^k$, say $x^{k_1}, \ldots ,x^{k_n}$. 
By results of Bianconi \cite{B05}, however, a function $x^{\beta}$ is definable in a structure only if $\beta$ is in the field generated by $k_1, \ldots ,k_n$. Since we can find such a $\beta$, and $r^x$ can be used to define $x^{\beta}$, we must conclude $r^x$ is not defined anywhere in $\calR$. 
Hence by the above proposition, the model companion exists.
\end{proof}
\end{ex}

\section{Neostability and Tameness}\label{stability}

In \cite{KTW18}, Kruckman, Tran, and Walsberg also give many general criteria under which neostability properties, such as NIP and NSOP$_1$, are preserved by the fusion of theories with these properties.
These preservation results are very much in the spirit of \cite{CP98}, \cite{E18}, and \cite{KR18}.
However, they also provide examples of interpolative fusions that have TP$_2$ despite both theories being fused having NTP$_2$.
In this section, we similarly show that depending on the o-minimal base theory $T$, the model companion can have IP or NIP, can be strong or not strong, and can have NTP$_2$ or TP$_2$.

For a comprehensive list of definitions and equivalent formulations of the properties NIP, NTP2, and strong or finite burden, please refer to \cite{S15} and \cite{A07}, respectively.
We first observe that if $T$ is just the theory of ordered divisible abelian groups and $D=M$, then the theory $T_{\calG}$ is a ``dense pair'' in the sense of \cite{vdD98}.
In \cite{vdD98}, van den Dries shows that the theory $T^d$ of a dense pair is complete if $T$ is complete. 
Since $T_{\calG}$ aligns with $T^d$ when $T$ is the theory ODAG, we conclude that the theory $T_{\calG}^*$ is axiomatized by $T^d$ as defined in \cite{vdD98}.
Moreover, by \cite{BDO11} this theory has NIP (it is not hard to show that $T_{\calG}^*$ has the property they call ``innocuous'' in \cite{BDO11}).
On the opposite end of the spectrum, the model companion $T^*_{\calG}$ has TP$_2$ if $T$ extends the theory of real closed fields, $\calG$ is a subgroup of the multiplicative group, and $T_{\calG}$ is companionable.

\begin{thm}[Theorem \ref{multTP2}]\label{TP2}
For $T$ the theory of a real closed field, the theory $T_{\calG}$, where $\calG$ is a subgroup of the multiplicative group, has a model companion $T^*_{\calG}$ that has TP$_2$.
\begin{proof}
It follows from Theorem \ref{mult} that the theory $T_{\calG}$ has a model companion.
Let $\calM \models T_{\calG}^*$ be an $\aleph_1$-saturated model, 
and let $(b_i)_{i \in \N}$ and $(c_{(i,j)})_{i,j \in \N}$ be countable sequences of elements in $M$,  
where every finite subset of $(c_{(i,j)})_{j \in \N}$ is $\Q$-linearly independent over $\calG$.
We can choose the sequence $(b_i)_{i \in \N}$ to be $\Q$-linearly independent from each other and the sequence $(c_{(i,j)})_{(i,j) \in \N^2}$ over $\calG$, since as a divisible subgroup $\calG$ must be infinite index subgroup of ${\bf D}$, and since we take $\calM$ to be suitably saturated.
For any $n \in \N$, let $A$ be the $n \times 2n+1$-matrix given as follows:
$$A_{i,k} = \begin{cases} b_i & k=1,\\ 
-1 & 1<k<n+2 \land k=i+1\\ 
1 & n+2\leq k\leq 2n+1 \land k=n+i+1\\  
0 & \text{otherwise}\end{cases}.$$ 

Consider the array generated by these indiscernible sequences and the formula $\vPhi(x, b_i,c_{(i,j)}) = b_i \cdot x + c_{(i,j)} \in \calG$.
Any two formulas in a row of the array are inconsistent, since each column requires that $b_i \cdot x$ is in a different coset of $\calG$.
By the $\Q$-independence of the sequences $(b_i)_{i \in \N}$ and $(c_{(i,j)})_{(i,j) \in\N^2}$ over $\calG$, 
for any $n \in \N$,
the corresponding matrix $A$ has rank $n$, which implies n+1 free variables for the corresponding solution set.
So given $j_1< \ldots <j_n \in \N$, there are infinitely many tuples $(x,y_1, \ldots ,y_n, c_{(1,j_1)}, \ldots ,c_{(n,j_n)})$ that satisfy the equation $A \cdot (x,y_1,\ldots,y_n, c_{(1,j_1)}, \ldots ,c_{(n,j_n)}) = \vec{0}$.
Rewriting the matrix so that the sub-tuple $(c_{(1,j_1)}, \ldots ,c_{(n,j_n)})$ is regarded as the part uniquely determined by the system of equations, we can then regard the sub-tuple $(x,y_1, \ldots,y_n)$ as the free variables.
In particular, we conclude that the set of components $(x,y_1, \ldots ,y_n)$ that correspond to a fixed $(c_{(1,j_1)}, \ldots ,c_{(n,j_n)})$ as a solution set for the above matrix equation has interior in $M^{n+1}$.
Hence the companion axioms tell us we can find a solution such that $y_1, \ldots ,y_m \in \calG$.
So $T^*_{\calG}$ has TP$_2$, as witnessed by this array.
\end{proof}
\end{thm}

Let $VS$ be the theory of an ordered vector space over $K$, a subfield of $\R$, in the language $\calL = \{<,0,+,(k)_{k \in K}, (k(x))_{k \in K} \}$, where $(k)_{k \in K}$ enumerates constant symbols for each element of $K$, and $(k(x))_{k \in K}$ enumerates scalar multiplication functions for each element of $K$.
In particular, $VS$ contains the axioms for an ordered divisible abelian group as well as the axioms for an ordered vector space.
In the language $\calL$ we know $VS$ has quantifier elimination and $VS_{\calG}$ does not have UEP, yet  $VS_{\calG}$ is not model complete by Example \ref{VSex}.
However, it does have a model companion $VS^*_{\calG}$. 
Using quantifier elimination we can show that $VS^*_{\calG}$ has NIP and may or may not have finite burden, depending on the base field.
To show both NIP and that finite burden occurs in a special case, we first need the following lemma.
Below, we use ``$\cdot$'' to denote the usual dot product.

\begin{defn}
Let $VS$ be the theory of an ordered vector space over $K\subseteq \R$, and let $(\calM, \calG) \models VS^*_{\calG}$. 
Define a \emph{coset type} of an element $a \in M$ over $C \subseteq M$ to be a maximal consistent set of formulas that are either of the following form or its negation:
$$ua-\vec{v}\cdot \vec{c} \in \calG $$
where $u,\vec{v} \in K$, $\vec{c} \in C$, and each of these formulas is true in $(\calM, \calG)$.
\end{defn}

\begin{lem} \label{quantred}
Let $VS$ be the theory of an ordered vector space over $K\subseteq \R$.
For every $(\calR, \calG) \models VS^*_{\calG}$ and for each $r \in R$, the following hold:
\begin{enumerate}
\item{For any $\acl$-independent set $C \subseteq R$,  the type $\tp (r/C)$ is implied by the $<$-cut of $r$ in $\acl (C) := K \langle C \rangle$ in conjunction with the coset type over $C$.
}
\item{Suppose $K = \Q( \eta_1, \ldots ,\eta_n)$ is a finite-dimensional extension of $\Q$ as a vector space.
Then for every $(\calR, \calG) \models VS^*_{\calG}$, for each $r \in R$ there is a finite set of elements $\{d_1, \ldots, d_n\} \subseteq R$ such that for any $C \subseteq R$ containing these elements,  the type $\tp (r/C)$ is implied by the $<$-cut of $r$ in $K \langle C \rangle$ plus the coset type over $\{d_1,\ldots ,d_n\}$.}
\end{enumerate}
\begin{proof}
For (1), let $(\calR, \calG) \models VS^*_{\calG}$ and $\{r\}, C \subseteq R$ be as in the hypotheses, and suppose $\varphi (x, \vec{c}) \in \tp (r/C)$.
By model completeness, this is equivalent to a disjunct of formulas of the form 
$$\exists \vec{y} \psi(x, \vec{y}, \vec{c}) \land \bigwedge_{i \in I} k_i x + \vec{u}_i \cdot \vec{y} + \vec{v}_i \cdot \vec{c} \in \calG \land \bigwedge_{i \not \in I} k_i x + \vec{u}_i \cdot \vec{y} + \vec{v}_i \cdot \vec{c} \not \in \calG $$
where $m \in \N$, $I \subseteq [m]$, for all $i \in [m]$ we have $k_i,\vec{u}_i,\vec{v}_i \in K$, and $\psi$ is a quantifier-free $\calL$-formula without disjuncts.
By quantifier elimination for ordered real vector spaces, we know that every definable function used in $\psi(x, \vec{y}, \vec{c})$ is a $K$-linear function.

Since $VS$ eliminates $\exists ^{\infty}$, either there is $\ell \in \N$ such that $\calR \models \exists^{\leq \ell} \vec{y} \psi(r, \vec{y}, \vec{c})$ or else $\calR \models \exists^{\infty}\vec{y} \psi(x,\vec{y},\vec{c})$.
If the former holds, we can write each $y_i$ as a $K$-linear function of $r$ and $\vec{c}$, and the subformula
\begin{gather} \label{lineq}
\bigwedge_{i \in I}  k_ir + \vec{u}_i \cdot \vec{y} + \vec{v}_i \cdot \vec{c} \in \calG \land  \bigwedge_{i \not \in I}  k_ir + \vec{u}_i \cdot \vec{y} + \vec{v}_i \cdot \vec{c} \not \in \calG 
\end{gather}
of $\varphi(r,\vec{c})$ is implied by $ \bigwedge_{i \in I} k'_i r + \vec{v}'_i \cdot \vec{c} \in \calG \land \bigwedge_{i \not \in I} k'_i r + \vec{v}'_i \cdot \vec{c} \not \in \calG$ for some other $k'_i\in K$ and $\vec{v}'_i \in R^{|\vec{c}|}$ for each $i \in \{1, \ldots ,m\}$.
Subtracting $\vec{v}_i \cdot \vec{c}$ on each side, we see that these conjuncts form a coset type over $C$, which we see implies $\tp (r/ C)$ in conjunction with the $<$-cut of $r$ in $K \langle C \rangle$.

Suppose now that $\calR \models \exists^{\infty} \vec{y} \psi(r,\vec{y},\vec{c})$.
Without loss of generality, we assume there is a witness for $\exists \vec{y} \psi(r,\vec{y},\vec{c})$ such that $\vec{y}$ is $\acl_{\calL}$-independent over $C \cup \{ r \}$.
Otherwise, we could write the $s< |\vec{y}|$ dependent coordinates as a $\{r\} \cup C$-definable function of the $s$ independent coordinates of $\vec{y}$, which can be subsumed into $\psi$.
Hence, by modifying the way we express $\vPhi (x,\vec{c})$ slightly to be in the appropriate form, we can find an  axiom in the model companion axiom scheme of $VS^*_{\calG}$ which
tells us that the formula $\exists \vec{y} \psi (x,\vec{y},\vec{c})$ implies $\varphi(x,\vec{c})$. 
So $\vPhi(x,\vec{c})$ is implied by an $\calL(C)$-formula, as desired.
Therefore this formula is already implied by part of the cut of $r$ over $K \langle C \rangle$, so we are done.

For (2), we show that in the special case that $K$ is finite-dimensional over $\Q$ as a vector space, we get the further quantifier reduction from a similar analysis of the formulas in $\tp (r/C)$.
Above in equation \ref{lineq}, each $k_i$, and each component of $\vec{u_i}$ and $\vec{v}_i$ is equal to $q_{i,0} + q_{i,1 }\eta_1 + \ldots +q_{i,n}\eta_n $ for some $q_{i,j} \in \Q$ for each $i \in \{1, \ldots ,m\}$ and $j \in \{0, \ldots ,n\}$, 
and for notational convenience define $\eta_0 = 1$.

Hence $ \bigwedge_{i \in I} k'_i r + \vec{v}'_i \cdot \vec{c} \in \calG \land \bigwedge_{i \not \in I} k'_i r + \vec{v}'_i \cdot \vec{c} \not \in \calG$
is implied by 
$$\bigwedge_{j \in I'} \eta_jr+\vec{v}'_j \cdot \vec{c}' \in \calG \land \bigwedge_{j \not \in I'} \eta_jr+\vec{v}'_j \cdot \vec{c}' \not \in \calG $$ 
for some $I' \subseteq \{0, \ldots ,n\}$.
Moreover, by the hypotheses each of the negated subformulas $\eta_i x+ \vec{v}'_i \cdot \vec{c} \not \in \calG$ is implied by any formula of the form $\eta_i x -d_{i+1} \in \calG$ for $i \in \{0, \ldots ,n\}$, one of which holds for $r$.
\end{proof}
\end{lem}

That $VS^*_{\calG}$ has NIP follows from an analysis of indiscernible sequences in light of this quantifier reduction.

\begin{thm}\label{VSNIP}
For $VS$ the theory of an ordered vector space over $K$ a subfield of $\R$, the theory $VS_{\calG}^*$ has NIP.
\begin{proof}
Let $(\calR, \calG) \models VS^*_{\calG}$ be a monster model, though we will only use that it is $|K|^{+}$-saturated.
We suppose for contradiction that there is a formula $\varphi(x, \vec{y})$ along with an element $a \in R$ and indiscernible sequence $( \vec{b}_i )_{i< \omega}$ that witnesses IP for $VS^*_{\calG}$,  i.e. $(\calR, \calG) \models \varphi(a, \vec{b}_i)$ precisely if $i$ is even. Let $|\vec{y}| = n$.
By model completeness of $VS^*_{\calG}$, the formula $\varphi$ is equal to a disjunct of formulas of the form
$$\sigma(x,\vec{y}) := \exists \vec{z} \Big ( \psi(x, \vec{y}, \vec{z}) \land \bigwedge_{j \in I} k_j x + \vec{u}_j \cdot \vec{y} + \vec{v}_j \cdot \vec{z} +c_j \in \calG \land \bigwedge_{j \not \in I} k_j x + \vec{u}_j \cdot \vec{y} +\vec{v}_j \cdot \vec{z}+c_j \not \in \calG \Big ) $$
where $I \subseteq [m]$ and $|\vec{z}| = d$ for some $m,d \in \N$, each $k_j, \vec{u}_j, \vec{v}_j$ and $c_j$ is in $K$, and $\psi$ is a quantifier free $\calL$-formula without disjuncts.
Since NIP is preserved under boolean combinations, one such disjunct must itself witness IP.
For convenience of notation, we will change the conjunct 
$$\bigwedge_{j \in I} k_j x + \vec{u}_j \cdot \vec{y} + \vec{v}_j \cdot \vec{z} +c_j \in \calG $$
to $A (x, \vec{y}, \vec{z}) + \vec{c} \in \calG^{|I|}$
 where $A$ is the matrix representation for the EIC transformation that corresponds to the concatenation of the linear transformations appearing in the specified conjuncts.
Similarly for $\bigwedge_{j \not \in I} k_j x + \vec{u}_j \cdot \vec{y} +\vec{v}_j \cdot \vec{z}+c_j \not \in \calG $ and $A' (x, \vec{y}, \vec{z}) + \vec{c}' \in (\calG^c)^{|I^c|}$.

Since $VS$ has NIP, either $\calR \models \exists^{\infty} \vec{z} \psi(a, \vec{b}_i, \vec{z})$ for cofinitely many $i< \omega$, or $\calR \models \exists^{\leq N} \vec{z} \psi(a, \vec{b}_i, \vec{z})$ for cofinitely many $i< \omega$.
We consider the first case.
Let $a \in M$ and $(b_i)_{i<\omega}$ be such that for cofinitely many $i < \omega$ there are infinitely many $\vec{r}$ in $R^{|\vec{z}|}$ that witness $\exists \vec{z} \psi (a,\vec{b}_i,\vec{z})$.
Without loss of generality, we suppose this holds for all $i< \omega$.
By cell decomposition, we can assume that the set of $\vec{z}$ that satisfy $\psi(a,\vec{b}_i,\vec{z})$ has interior in $R^d$ for cofinitely many $i<\omega$.

We now adjust the form that the formula $\sigma(x, \vec{y})$ takes in order to apply the model companion axioms as listed in section \ref{(R,G)}.
We can replace the tuple $\vec{z}$ with the tuple $\tilde{z} = (\vec{z}, \vec{z}')$ where $\ell:=|\vec{z}'| = |I|$.
We consider the following formula that is equivalent to $\sigma(x,\vec{y})$:
$$ \tilde{\sigma}(x,\vec{y}) := \exists \tilde{z} \psi(x,\vec{y},\tilde{z}_{[d]})\land \Big( A \tilde{z}_{[d]} +\vec{c} = (\tilde{z}_{d+1}, \ldots ,\tilde{z}_{d +\ell})\Big) \land \bigwedge_{j=d+1}^{d+\ell} \tilde{z}_j \in \calG
 \land\Big( A' \tilde{z}_{[d]} + \vec{c}' \not \in \calG^{|I^c|} \Big)$$
where $\tilde{z}_{[d]} = (\tilde{z}_1, \ldots ,\tilde{z}_d)$.
Now we can apply the model companion axiom in which the $\calL$-formula is $\psi(x,\vec{y},\tilde{z}_{[d]})\land A \tilde{z}_{[d]} +\vec{c} = (\tilde{z}_{d+1}, \ldots ,\tilde{z}_{d +\ell})$ and $\vec{x}_0 = (\tilde{z}_{d+1}, \ldots ,\tilde{z}_{d+\ell})$, and the functions $(\tau_i)_{i<\ell}$ are given by $A' \tilde{z}_{[d]} + \vec{c}'$.
The corresponding model companion axiom tells us that $\tilde{\sigma}$ holds precisely if $\theta_{\psi}$ does, which then therefore holds precisely if $\sigma$ does.
Since $\theta_{\psi}$ is purely an $\calL$-formula, NIP for $VS$ tells us that
$\theta_{\psi}$ holds for cofinitely many $i< \omega$, or does not hold for cofinitely many $i<\omega$, but this contradicts $\varphi$ having IP.

We now consider case two, that for some $N \in \N$ we have $\calR \models \exists^{\leq N} \vec{z} \psi (a, \vec{b}_i, \vec{z})$ for cofinitely many $i< \omega$.
Without loss of generality, we suppose this is true for all $i < \omega$.
If there are at most $N$ elements of $R^{|\vec{z}|}$ that witness $\calR \models \exists \vec{z} \psi (a, \vec{b}_i, \vec{z})$, then each such witness is in $\dcl_{\calL}(\{a\} \cup \vec{b}_i)$.
So by o-minimality we can enumerate them as $\calL$-definable functions of $a$ and $\vec{b}_i$ and whichever parameters appear in $\psi$, say as $\vec{f}_1(a, \vec{b}_i), \ldots ,\vec{f}_N(a, \vec{b}_i)$.
We conclude that for $x=a$ and $\vec{y} = \vec{b}_i$ for any $i< \omega$ the following sentence:
$$\bigvee_{\ell=1}^N \Big ( \psi(a, \vec{b}_i, \vec{f}_{\ell}(a,\vec{b}_i)) \land \bigwedge_{j \in I} k_j a + \vec{u}_j \cdot \vec{b}_i + \vec{v}_j \cdot \vec{f}_{\ell}(a,\vec{b}_i) \in \calG \land \bigwedge_{j \not \in I} k_j a + \vec{u}_j \cdot \vec{b} + \vec{v}_j \cdot \vec{f}_{\ell}(a,\vec{b}_i) \not \in \calG \Big ) $$
holds precisely if $\varphi(a, \vec{b}_i)$ does.

By the pigeonhole principle and by restricting (if necessary) to some cofinal subset of $\omega$, we may assume that the only part of $\varphi(a, \vec{y})$ that alternates in truth value on a cofinal subset of $\omega$ is some subformula corresponding to $\ell \in [N]$ of the form
$$ \bigwedge_{j \in I} k_j a + \vec{u}_j \cdot \vec{y} + \vec{v}_j \cdot \vec{f}_{\ell}(a,\vec{y}) \in \calG \land  \bigwedge_{j \not \in I} k_j a + \vec{u}_j \cdot \vec{y} + \vec{v}_j \cdot \vec{f}_{\ell}(a,\vec{y}) \not \in \calG.$$
By pigeonhole principle, at least one of the conjuncts in this subformula alternates in truth value on a cofinal subset of $\omega$, and we may make the further assumption that one such conjunct holds precisely if $i< \omega$ is odd.
Without loss of generality, suppose that $k_1x+ \vec{u}_1 \cdot \vec{y} + \vec{v}_1 \cdot \vec{f}_1(x,\vec{y}) \in \calG$ is such a subformula.
Since we know that $f_1(x,\vec{y})$ is a definable function in the vector space language, by quantifier elimination for ordered vector spaces we know that $f_1(x,\vec{y}) = u'x+ \vec{v}'\cdot \vec{y}+\vec{c}_1$ for some $u',\vec{v}' \in K$.

Hence the subformula is equal to $k_1x+ \vec{u}_1 \cdot \vec{y} + \vec{v}_1 \cdot (u'x+ \vec{v}' \cdot \vec{y}+\vec{c}_1)\in \calG$, which is furthermore equivalent to $k^*x + \vec{v}^{*} \cdot \vec{y} + c^* \in \calG$ for the requisite $k^*,\vec{v}^{*},c^*$ in $K$.
Hence if $i<\omega$ is odd, then
$\vec{v}^{*} \cdot \vec{b}_i  \in \calG - (k^*a + c^*)$.
However, if $i<\omega$ is even, then there exists $c_i \in (\calG + c^*)^c$ such that
$\vec{v}^{*} \cdot \vec{b}_i  \in \calG - (k^*a + c_i)$.
So if $i<\omega$ is even, then $\vec{v}^{*} (\vec{b}_i - \vec{b}_{i+1}) \in \calG - (c_i - c^*)$, where $c_i - c^* \not \in \calG$.
However, $\vec{v}^{*} (\vec{b}_{i+1} - \vec{b}_{i+3}) \in \calG $ since $\vec{v}^{*}b_{i+1}$ and $\vec{v}^{*}b_{i+3}$ are in the same coset by the formula that the odd index $\vec{b}_j$'s satisfy.
By the indiscernibility of the sequence $(b_i)_{i< \omega}$, for any $j>i$ we must thereby conclude that $\vec{v}^{*} (\vec{b}_i - \vec{b}_{j}) \in \calG$, a contradiction.
So no such subformula can alternate in truth value on such an indiscernible sequence.
By preservation of NIP under boolean combinations of formulas, we conclude that $\sigma$ and hence $\vPhi$ has NIP, as desired.

\end{proof}
\end{thm}

Recall that the notion of ``finite burden'' as defined in \cite{A07} is equivalent to finite inp-rank.

\begin{cor}\label{VSn}
For $VS_n$, the theory of an ordered vector space over $K = \Q \langle \eta_1, \ldots ,\eta_n \rangle \subset$, i.e. $K$ is an algebraic extension with linear degree $n$ over $\Q$, the model companion $VS^*_{n,\calG}$ of $VS_{n,\calG}$ has burden at most $n+2$.
\begin{proof}
By the Lemma \ref{quantred}, for any $\calM \models VS^*_{n,\calG}$ and for any $x \in \calM$ and $C \subseteq M$ countable, the type of $x$ over $C$ is determined by the $\calL$-type of $x$ over $C$ in conjunction with the formulas $x \in \calG + y_1$ and $\eta_1x \in \calG + y_2, \ldots , \eta_nx \in \calG + y_n$ for some $y_1, \ldots ,y_n$ in $K \langle C \rangle$.
We suppose for contradiction that there is an inp-pattern of depth $n+3$, and that the array of formulas $\langle \psi_i(x,y): i < n+3 \rangle$ and the array of indiscernible sequences $(c_{i,j})_{i \in [n+3],j \in \N}$ witnesses this.
In particular, assume that the first $n+2$ rows form an inp-pattern of depth $n+2$ on their own.

By dp-minimality of o-minimal structures, we know that the purely $\calL$-definable part of $\psi_i(x,y)$ is trivial (or is exactly the same for all $i \in [n+2]$ and is independent of the parameter tuple $y$) for all but one $i \in [n+2]$.
Without loss of generality let $\psi_0(x,y)$ be the lone formula in the array with a nontrivial $\calL$-definable component that varies as $y$ does.
The rest of the formulas must be definable in $\calL_{\calG}$ but not in $\calL$, and without loss of generality we may assume they do not define intervals of any kind.

By the quantifier reduction result Lemma \ref{quantred}, we know that each $\psi_i$ is equivalent to an $\calL$-formula or the disjunct of conjuncts of $\calL$-formulas with some conjuncts of the form $\eta_{\ell} x + \vec{v} \cdot y + k \in \calG$ or the form $\eta_{\ell} x + \vec{v} \cdot y + k \not \in \calG$, where $\vec{v}$ and $k$ are in $K$.
Since for $i \neq 0$ the purely $\calL$-definable part of $\psi$ is trivial, we conclude that $\psi_i(x,y)$ defines a finite union of a finite intersection of cosets of $\calG$ for one of $x, \eta_1x, \ldots,\eta_n x$ to lie in, and coset-complements for $x$ or one of those scalar multiples not to lie in.
Since $\calG$ is an infinite-degree subgroup of $\calM$, if $\psi_i$ were purely a finite union of finite intersections of coset-complements for $x$ and its scalar images, then for any $j_1, \ldots ,j_{\ell} <\omega$ with $\ell \in \N$ and $j_1< \ldots < j_{\ell}$ we could find an $x$ that satisfies all of $\psi_i(x, c_{i,j_1}), \ldots ,\psi_i(x, c_{i,j_{\ell}})$.
So to have finite inconsistency across each row of the array, for all but finitely many $i< \omega$ each disjunct of $\psi_i(x,y)$ must include a conjunct that for any given $y$ dictates the coset of $\calG$ for $x$ or a scalar multiple of $X$, and that coset varies as $y$ does.

We know from part (2) of Lemma \ref{quantred} that for any element $r \in M$ there are at most $n+1$ disjoint formulas needed to define the cosets of $\calG$ in which $r,\eta_1r,\ldots ,\eta_nr$ each lie respectively.
In conjunction with the $<$-cut and formulas excluding $r$ or some $\eta_i r$ from being in any $C$-definable coset of $\calG$, this isolates the type of $r$ over $M$.
Suppose that for each $i>0$, each formula $\psi_i$ defines a disjunct of cosets of $\calG$ for some subset of $\{x,\eta_1x, \ldots,\eta_nx\}$ to lie in.
In order that the paths of the inp-pattern be consistent, each $\psi_i$ can define the coset of at most one element of $\{x,\eta_1x, \ldots,\eta_nx\}$.
If $\Psi$ is a set of formulas such that each one defines the coset of at most one of $x,\eta_1x, \ldots,\eta_nx$, then $\Psi$ ony contains at most $n+1$ mutually non-equivalent formulas.
So each $\psi_i$ is equivent to one in $\Psi$, a set of at most $n+1$ distinct formulas.
Hence if $\psi_{n+1}(x, y)$ is the formula for the $n+2^{th}$ row, where $y=c_{(n+2),j}$ is the parameter used for the $j^{th}$ column, then by Lemma \ref{quantred} the formula $\psi(x, c_{n+2,j})$ is implied by formulas from the previous rows.
Thus the values of the parameters $c_{n+2,j}$ are dictated by formulas and parameters of the previous $n+1$ rows, contradicting indiscernibility of $(c_{i,j})_{i,j \in \N}$.
\end{proof}
\end{cor}

Below we see that this $VS^*_{\calG}$ is not strong in the sense of \cite{A07} if the base field $K$ has infinite linear degree over $\Q$.
The proof that $VS^*_{\calG}$ is it not strong is directly analogous to the above proof that the model companion for an expansion of a real closed field has TP$_2$.

\begin{rmk}\label{VSinf}
For $VS_{\infty}$ the theory of an ordered vector space over $K \subseteq \R$, the theory $(VS_{\infty})_{\calG}$ has a model companion $(VS_{\infty})^*_{\calG}$, and $(VS_{\infty})^*_{\calG}$ is not strong.
\begin{proof}
Consider the array of formulas where the formula with coordinates $(i,j)$ is $k_ix + c_{(i,j)} \in \calG$, where $k_i \in K$ and $c_{(i,j)} \not \in \calG$.
For each $i \in \N$, let $(c_{(i,j)})_{j \in \N}$ be a sequence of constants such that for all $k>j$ it is not the case that 
$c_{(i,k)}-c_{(i,j)} \in \calG$.
This can be arranged by the fact that $\calG$ has infinite $\Q$-linear degree over $K$, by divisibility and saturation.
Thus any two formulas in a row of the array are inconsistent.
Since $K$ has infinite linear degree over $\Q$, we can arrange that the set $A$ of elements from $K$ that appear in these formulas are $\Q$-linearly independent.
By saturation, we can arrange that the array $(c_{(i,j)})_{(i,j) \in \N^2}$ is such that the tuples $\{(k_i, c_{(i,j)}): i,j \in \N \}$ are $K$-linearly independent as well.
Hence it follows that for any $m \in \N$ and any $j_1, \ldots ,j_m \in \N$ it is true that the set of equalities $\{ k_1(x) + c_{(1,j_1)} = y_1, \ldots ,k_m(x) + c_{(m,j_m)} = y_m \}$ has infinitely many solutions for $(x,y_1, \ldots ,y_m)$, and the solution space as we vary the elements $c_{(i,j_i)}$ has linear degree $2m+1$ over $\Q$.
Hence we can rewrite the conjunct of these formulas to apply the companion axioms, which tell us for each $(c_{(1,j_1)}, \ldots ,c_{(m,j_m)})$ we can find a solution such that $y_1, \ldots ,y_m \in \calG$.
So $(VS_{\infty})^*_{\calG}$ is not strong, as witnessed by this array.
\end{proof}
\end{rmk}

Note that in an analogous manner to the proof of Remark \ref{VSinf}, we can always construct an inp-pattern of depth $n+2$ for the theory $VS^*_{n,\calG}$.
Hence Theorem \ref{VSNIP}, Corollary \ref{VSn}, and Remark \ref{VSinf} together constitute Theorem \ref{strongVS}.
As we have assumed that $T$ is o-minimal and as the theory $T^*_{\calG}$ is a model-complete extension of this theory, it is noteworthy that $T^*_{\calG}$ does not necessarily exhibit model-theoretic tameness.
The lack of correlation of the model companion $T^*_{\calG}$ with any of the widely employed neostability properties suggests a need for a more robust notion of tameness that captures the kind that the theory $T^*_{\calG}$ exhibits as an extension of $T_{\calG}$.

\subsection*{Acknowledgements} The author would like to thank Erik Walsberg for offering the statement of the main theorem as an apt conjecture, and for many conversations about the results which yield the proof of said conjecture. 
Many thanks go to Philipp Hieronymi, who helped fix a number of the arguments and gave diligent attention to correcting many drafts of this paper. 
Great thanks also to Christian d'Elb\'ee for a detailed and helpful discussion of the main theorem, and for suggesting a nice proof of Theorem \ref{TP2}.
The author would also like to thank the anonymous referee for their attentive corrections and suggestions.
Finally, thanks to Chieu-Minh Tran and Elliot Kaplan for some very helpful discussions of various results in this paper.
This material is based upon work supported by the National Science Foundation Graduate Research Fellowship Program under Grant No. DGE -- 1746047.

%\nocite{*}

%%%%%%%%%%%%%%%%%%%%%%%%%%%%%%%%%%%%%%%%%%%%%%%%%%%%%%%%%%%%%%%%%%%%

\end{document}